# Algebraic estimates, stability of local zeta functions, and uniform estimates for distribution functions

By D. H. Phong and Jacob Sturm *

## Abstract

A method of "algebraic estimates" is developed, and used to study the stability properties of integrals of the form $\int_B |f(z)|^{-\delta} dV$, under small deformations of the function $f$. The estimates are described in terms of a stratification of the space of functions $\{R(z) = |P(z)|^{\varepsilon}/|Q(z)|^{\delta}\}$ by algebraic varieties, on each of which the size of the integral of $R(z)$ is given by an explicit algebraic expression. The method gives an independent proof of a result on stability of Tian in 2 dimensions, as well as a partial extension of this result to 3 dimensions. In arbitrary dimensions, combined with a key lemma of Siu, it establishes the continuity of the mapping $c \to \int_B |f(z,c)|^{-\delta} dV_1 \cdots dV_n$ when $f(z,c)$ is a holomorphic function of $(z,c)$. In particular the leading pole is semicontinuous in $f$, strengthening also an earlier result of Lichtin.

## 1. Introduction

The main purpose of this paper is to study integrals of the form

$$(1.1) \qquad \int_{B^{(n)}} |f(z)|^{-\delta} dV_1 \cdots dV_n,$$

where $B^{(n)}$ is a polydisk in $\mathbb{C}^n$, $f(z) = (f_j(z))_{j=1}^J$ is a $J$-dimensional vector of holomorphic functions on $B^{(n)}$, and $|f(z)|^2 = \sum_{j=1}^J |f_j(z)|^2$. Such integrals are sometimes referred to in the literature as "local zeta functions." They have long been a subject of investigation in various branches of mathematics. Their basic analytic properties have been established in the work of Bernstein and Gel'fand [4], Atiyah [3], Gel'fand and Shilov [6], Arnold, Gussein-Zadé, Varchenko [2], and Igusa [8], among others. As is the case with the Riemann zeta function, the integral converges for $\delta$ in a half-plane. It has a meromorphic continuation to the entire plane, with the location of its poles easily read off from the resolution

*Research supported in part by the National Science Foundation under NSF grant DMS-98-00783.

of singularities of $f(z)$. In singularity theory [2], [15], [17], the inverse of the supremum $\delta_0$ over all exponents $\delta$ for which the integral (1.1) is finite provides a natural notion of multiplicity for $f(z)$. In harmonic analysis, the integral (1.1) controls the distribution function of $f(z)$ via the Chebychev inequality:

$$(1.2) \qquad \mathrm{Vol}\{z \in B^{(n)}; \; |f(z)| \le \alpha\} \le \alpha^\delta \int_{B^{(n)}} |f(z)|^{-\delta} dV,$$

and bounds on integrals of type (1.1) have recently emerged as central to aspects of complex differential geometry [1], [19], in particular to the solution of certain complex Monge-Ampère equations [20], [22], [24].

We shall be mainly concerned with the issue of stability for the finiteness of integrals of type (1.1), under deformations of the function $f(z)$. Surprisingly little is known, and the only results available to date appear to be the following: Tian [22] has shown that, in two dimensions, the finiteness of the integrals (1.1) is stable under holomorphic perturbations of $f(z)$ with small sup norms. As part of an unpublished work on Kähler-Einstein metrics for Fano manifolds, Siu [20] has shown that, in arbitrary dimensions, if (1.1) is finite for $f(z) = f(z, 0)$ where $f(z, c)$ holomorphic in $(z, c) \in B^{(n)} \times B^{(1)}$, then the infimum over $\{0 < |c| < \rho\}$ of such integrals for $|f(z, c)|^{-\delta}$ is finite for all $\rho$ sufficiently small. A closely related version of this lemma of Siu has appeared in his work on the Fujita conjecture [19][1]. Lichtin [14] has shown that under conditions similar to Siu's, the integrals (1.1) remain finite for all $|c|$ small enough, but with the additional assumption that $f(z, c)$ have an isolated singularity for all $c$, and $B^{(n)}$ be replaced by a Milnor ball $B^{(n)}(c)$ which may be $c$-dependent.

The case of $f(z)$ real-analytic and $B^{(n)} \subset \mathbb{R}^n$ is somewhat better understood, but it has been unclear whether it is any reliable guide for the holomorphic case. In the real case, it is known that (1.1) is stable in 2 real dimensions, thanks to a theorem of Karpushkin [11]–[12], but not in dimensions 3 or higher, thanks to the following counterexample of Varchenko [23]

$$(1.3) \qquad f(x, \varepsilon) = (x_1^4 + \varepsilon x_1^2 + x_2^2 + x_3^2)^2 + x_1^{4p} + x_2^{4p} + x_3^{4p}.$$

The finiteness of (1.1) is unstable for $f(x, \varepsilon)$, since $\delta_0 = \frac{5}{8}$ for $\varepsilon = 0$, $\delta_0 = \frac{3}{4}$ for $\varepsilon > 0$, and $\delta_0 < \frac{1}{2} + \gamma(p)$ for $\varepsilon < 0$ and $\lim_{p \to \infty} \gamma(p) = 0$. Note that this example does not rule out the possibility of a result analogous to Siu's for the real case.

In this paper, we develop a new method for the study of stability of integrals of holomorphic functions. A key component of the method is certain uniform estimates and stability for complex integrals of "rational" expressions of the form

$$(1.4) \qquad R(z) = \frac{(\sum_{i=1}^I |P_i(z)|^2)^{\varepsilon/2}}{(\sum_{j=1}^J |Q_j(z)|^2)^{\delta/2}} = \frac{|P(z)|^\varepsilon}{|Q(z)|^\delta},$$

where both $P_i(z)$ and $Q_j(z)$ are one variable polynomials of bounded degrees, and the domain of integration is a ball in the complex plane. It is well-known that integrals such as (1.4) are highly transcendental: they cannot be evaluated in closed form in general, and even in the exceptional cases where they can, they produce transcendental objects such as logarithms, inverse trigonometric, and elliptic functions. It may therefore be surprising that, paradoxically, *algebraicity is restored* if we focus not on the *exact* values of the integrals (1.4), but only on their *sizes*. In this sense, the uniform estimates we provide are "algebraic estimates." The fundamental fact which emerges is that the space of $(P(z), Q(z)) = (P_i(z), Q_j(z))$ can be stratified by constructible algebraic varieties, on each of which the size of the integrals of $R(z)$ can be expressed again by expressions of the form (1.4), but whose variables are now the coefficients of $P_i(z)$ and $Q_j(z)$, and whose coefficients are integers which do not depend on $R$ (see Theorem 4 for a precise statement). Some basic techniques for the method such as cluster scales were introduced for the study of the real case in [18]. There the method was used to give an independent proof of Karpushkin's theorem as well as a sharp stability theorem in 3 real dimensions which fits the constraints dictated by Varchenko's example. But the complex setting is the natural setting for the rational expressions (1.4) and their stratification by complex varieties, and it is here that algebraic estimates can be formulated in their full generality and that their underlying geometry becomes apparent.

We describe now our main results. In 2 dimensions, we obtain a new proof of Tian's result. In 3 dimensions, we obtain a new stability theorem under arbitrary holomorphic deformations, if the exponent $\delta$ in (1.1) satisfies the condition $\delta < 4/N$, where $N$ is the order of vanishing of $f(z)$ at the origin.

For arbitrary dimensions, making essential use of Siu's lemma and resolution of singularities, we obtain the following theorem which may be termed holomorphic stability for 1-parameter deformations:

MAIN THEOREM. *Let $g(z, c)$ be a $J$-vector of holomorphic functions on a polydisk $B^{(n)} \times B^{(1)}$, and assume that $\int_{B^{(n)}} |g(z, 0)|^{-\delta} dV_1 \cdots dV_n < \infty$. Then there exists a smaller polydisk $B'^{(n)} \times B'^{(1)}$ so that the function $c \to \int_{B'^{(n)}} |g(z, c)|^{-\delta} dV_1 \cdots dV_n$ is finite and continuous for $c \in B'^{(1)}$.*

The main theorem implies Lichtin's theorem and provides a strengthened version of Siu's lemma. It also shows that stability properties in the real setting vary sharply from those in the complex setting, since it rules out a complex version of Varchenko's counterexample.

Given the diversity of methods required in the works of Karpushkin, Lichtin, Siu, Tian, and Varchenko (which range from the versal theory of deformations to Carleman estimates for the $\bar{\partial}$ operator), it is encouraging that the present method has made contact with them all.

This paper is divided in two main parts, with the algebraic estimates developed in the first part, consisting of Sections 1–4, and the applications to stability and distribution functions developed in the second part, consisting of Sections 5–8. More precisely, in Section 2, we present estimates for the integrals of $|P(z)|^{\varepsilon}/|Q(z)|^{\delta}$ in terms of local cluster scales $L_k(\alpha)$ for the roots of $Q(z)$. Section 3 is devoted to the special case where $P(z)$ is a constant. The goal here is to "symmetrize" the cluster scale estimates, that is, to re-express the estimates in terms of rational expressions in the coefficients of $Q(z)$. Section 4 is devoted to the symmetrization problem for general rational expressions of type (1.4). Three useful techniques are introduced: The first is a regularization process, which allows us in effect to disentangle the zeroes of $P(z)$ and $Q(z)$. The second is the use of $\theta$-parameters, which allows us to replace a $J$-vector $f(z) = (f_j(z))_{j=1}^{J}$ in the denominator of the integrand by a scalar function $\sum_{j=1}^{J} e^{2\pi i \theta_j} f_j(z)$, at the expense of introducing a new integral over the $\theta$ domain. The third technique consists of sampling lemmas which reduce our considerations to a finite number of $\theta$ values. The key stratification $\mathcal{U}_{\lambda}$ is also described there, and the main result is presented in Theorem 4. In Section 5, we collect some general facts and definitions about stability. For our purposes, we require an extension of an important earlier result of Stein [21], which we establish using Hironaka's theorem on resolution of singularities. Section 6 is devoted to stability in dimensions $n \leq 3$. A characteristic feature of these results is that they only involve rational expressions of the form (1.4) with constant numerator. The Main Theorem is proved with the help of Siu's lemma in Section 7, in a more general form using plurisubharmonic functions (which also appear in Siu's work). In Section 8, we have listed some immediate consequences of our work for the stability of bounds for distribution functions. This is a topic of particular current interest, with some recent advances described in [5] and [18].

Finally, in the Appendix, we have reproduced with Professor Y.-T. Siu's kind permission the statement and proof of his unpublished result on holomorphic stability in arbitrary dimensions, which plays an essential role in this work.

*Acknowledgements.* The authors would like to thank E. M. Stein for earlier collaboration on related topics which led to this work, R. Friedman, M. Kuranishi, M. Robinson, and M. Thaddeus, for clarifications of different aspects of Hironaka's theorem, Y.-T. Siu and G. Tian for many stimulating conversations on stability and its applications to geometry. The authors are particularly grateful to Y.-T. Siu for informing them of his unpublished result on stability several years ago, and for allowing them to incorporate in this paper his proof of this result, which is crucial to their work.

## 2. Local cluster-scale estimates

The key ingredient of our approach is certain uniform estimates for integrals of rational functions. These estimates can be formulated either in terms of "local cluster scales" $L_i(\alpha)$ at each root $\alpha$ of the denominator, or in terms of analytic functions of the coefficients of both numerator and denominator. In this section, we derive the estimates in terms of cluster scales. First, we require some notation:

- Let $dV = dx\,dy$, the standard Euclidean measure on $\mathbb{C}$.

- For $r > 0$, let $B_r$ be the open disk of radius $r$ in $\mathbb{C}$, centered at 0. Sometimes we shall just write $B$ for the disk of radius 1.

- Let $P(z)$, $Q(z)$ be polynomials with complex coefficients of degrees $M$ and $N$ respectively, with $Q(z)$ monic. Let $S = \{\alpha : Q(\alpha) = 0\}$ denote the set roots of $Q(z)$, counted with multiplicity (so $S$ is a set with $N$ elements).

- If $A \subseteq \mathbb{C}$, define the diameter $d(A)$ of $A$ by $d(A) = \sup_{\alpha, \beta \in A} |\alpha - \beta|$.

- For $0 \leq k \leq N - 1$ and $\alpha \in S$ define

$$(2.1) \qquad L_k(\alpha) = \inf\{d(S_{N-k}(\alpha))\},$$

where the infimum is taken over all subsets $S_{N-k}(\alpha) \subseteq S$ such that $|S_{N-k}(\alpha)| = N - k$ and $\alpha \in S_{N-k}(\alpha)$. Observe that

$$(2.2) \qquad L_0(\alpha) \geq L_1(\alpha) \geq \cdots \geq L_{N-2}(\alpha) \geq L_{N-1}(\alpha) = 0,$$

and that $\alpha$ is a root of multiplicity $N - k$ if and only if $L_k(\alpha) = 0$.

Our goal is to estimate integrals of the form

$$(2.3) \qquad \int_{B_\Lambda} \frac{|P(z)|^\varepsilon}{|Q(z)|^\delta}\,dV$$

where $\varepsilon, \delta$ are nonnegative real numbers, and $\Lambda > 0$.

We shall make the following two assumptions:

$(2.4)$     (1) $S \subseteq B_{\Lambda/2}$.

         (2) $\nu\varepsilon + 2 - (N - k)\delta \neq 0$, for all integers $k, \nu$

             with $0 \leq k \leq N$, $0 \leq \nu \leq M$.

Assumption (2), which excludes finitely many lines in the $(\varepsilon, \delta)$ plane, is made in order to simplify the final form of our answer; it may be easily removed at the expense of introducing certain *log* terms in our estimates. However, our main applications require only the consideration of a dense set of rational values for $\varepsilon$ and $\delta$, so we shall omit these technicalities.

For each $\nu \geq 0$, we define $k_\nu$ by $k_\nu = -1$ if $\nu\varepsilon + 2 > N\delta$, and otherwise as the integer between $0$ and $N - 1$ satisfying

$$(2.5) \qquad (N - k_\nu - 1)\delta < \nu\varepsilon + 2 < (N - k_\nu)\delta.$$

Evidently, the integral (2.3) is finite if and only if $L_{k_\nu}(\alpha) > 0$ for any root $\alpha$ of $Q(\alpha) = 0$, where $\nu$ is the order of vanishing of $P(z)$ at $\alpha$, and $k_\nu$ is defined as above. Since the cluster scales $L_k(\alpha)$ are decreasing in $k$, this condition is actually equivalent to the seemingly more restrictive condition that for all $\nu$ with $P^{(\nu)}(\alpha) \neq 0$, and $k_\nu$ defined by (2.5), we have $L_{k_\nu}(\alpha) > 0$. The following theorem gives a precise, quantitative version of this statement.

THEOREM 1. *Under the preceding assumptions, the integral* (2.3) *is of size*

$$(2.6) \qquad \int_{B_\Lambda} \frac{|P(z)|^\varepsilon}{|Q(z)|^\delta}\, dV \sim \sum_{\{\alpha : Q(\alpha) = 0\}} \sum_{\{\nu : P^{(\nu)}(\alpha) \neq 0\}} \frac{|P^{(\nu)}(\alpha)|^\varepsilon}{\Phi_{\nu, k_\nu}(\alpha)}\,,$$

*where* $\Phi_{\nu, k}(\alpha)$ *is defined by*

$$(2.7) \qquad \Phi_{\nu, k}(\alpha) = \begin{cases} L_k(\alpha)^{(N-k)\delta - (\nu\varepsilon + 2)} \prod_{0 \leq i < k} L_i(\alpha)^\delta, & \text{if } k \geq 0; \\[2mm] \Lambda^{N\delta - (\nu\varepsilon + 2)}, & \text{if } k < 0, \end{cases}$$

*and for each* $\nu \geq 0$, $k_\nu$ *is defined as in* (2.5).

Here the equivalence $\sim$ means that each side is bounded by positive constant multiples of the other side, with constants which depend only on $\varepsilon$, $\delta$, and the degrees $M$ and $N$ of $P(z)$ and $Q(z)$. The constants are independent of the choice of $P$ and $Q$.

*Proof.* To prove the theorem we decompose the domain of integration $B_\Lambda$ as

$$(2.8) \qquad B_\Lambda = \cup_{\alpha \in S} D(\alpha)$$

where $D(\alpha)$ is defined by

$$(2.9) \qquad D(\alpha) = \{z \in B_\Lambda : \quad |z - \alpha| \leq |z - \beta| \quad \text{for all } \beta \in S\}.$$

We note that any $z \in B_\Lambda$ must be in $D(\alpha)$ for some $\alpha$. Furthermore, when $\alpha \neq \beta$, the intersection of $D(\alpha)$ and $D(\beta)$ is contained in a line, and is hence of measure 0. Thus the integral over $B_\Lambda$ may be written as a sum, over all $\alpha$, of the corresponding integrals over the $D(\alpha)$.

*Upper bounds.* We start by showing that the left side of (2.6) is less than or equal to a constant times the right side. To do this, we shall show that for each fixed $\alpha \in S$ the integral over $D(\alpha)$ is bounded above by a constant times the summand corresponding to $\alpha$ on the right side. Thus we fix a root $\alpha \in S$. We construct an ordering of the other roots of $Q(z)$ in the following manner: Choose $\beta_0 \in S$ such that $|\beta_0 - \alpha| \geq |\beta - \alpha|$ for all $\beta \in S$. For $1 \leq i \leq N-1$, choose $\beta_i \in S \backslash \{\beta_0, \cdots, \beta_{i-1}\}$ such that $|\beta_i - \alpha| \geq |\beta - \alpha|$ for all $\beta \in S \backslash \{\beta_0, \ldots, \beta_{i-1}\}$. We claim that for $0 \leq i \leq N-1$,

$$(2.10) \qquad L_i(\alpha) \sim |\beta_i - \alpha|,$$

with implied constants which depend only on $N$. We establish this estimate by induction in $i$. Using the definition of $L_0(\alpha)$, it follows immediately that $L_0(\alpha) \geq |\beta_0 - \alpha|$. On the other hand, we must have $|\beta_0 - \alpha| \geq L_0(\alpha)/2$. Otherwise $|\beta - \beta'| \leq |\beta - \alpha| + |\beta' - \alpha| \leq 2|\beta_0 - \alpha| < L_0(\alpha)$ for all $\beta, \beta' \in S$ which contradicts the definition of $L_0(\alpha)$. This establishes (2.10) for $i = 0$. Next assume (2.10) for all $j \leq i-1$. Since $S \backslash \{\beta_0, \cdots, \beta_{i-1}\}$ contains $\alpha$, we have $d(S \backslash \{\beta_0, \cdots, \beta_{i-1}\}) \geq L_i(\alpha)$. This implies that $|\beta_i - \alpha| \geq L_i(\alpha)/2$ (otherwise, just as we argued above, for any $\beta, \beta' \in S \backslash \{\beta_0, \cdots, \beta_{i-1}\}$, we would have

$$|\beta - \beta'| \leq |\beta - \alpha| + |\beta - \alpha'| \leq 2|\beta_i - \alpha| < L_i(\alpha),$$

which contradicts the definition of $L_i(\alpha)$.) To get the reverse inequality, we consider two cases: if $d(S \backslash \{\beta_0, \cdots, \beta_{i-1}\}) = L_i(\alpha)$, then $|\beta_i - \alpha| \leq L_i(\alpha)$ and we are done. If not, choose a subset $S_{N-i}(\alpha)$ of $N-i$ roots which achieves the minimum in (2.1). Then $L_i(\alpha) = d(S_{N-i}(\alpha))$ and there exists a $j < i$ such that $\beta_j \in S_{N-i}(\alpha)$. By induction, we have the estimates $L_i(\alpha) \leq L_j(\alpha) \sim |\beta_j - \alpha| \leq L_i(\alpha)$ and thus $L_j(\alpha) \sim L_i(\alpha)$. On the other hand, $|\beta_i - \alpha| \leq |\beta_j - \alpha|$. This completes the inductive step, and (2.10) is proved.

The estimate (2.10) leads to the following basic estimate in the region $D(\alpha)$ for each factor $|z - \beta_i|$ in the polynomial $Q(z) = \prod_{i=1}^{N}(z - \beta_i)$:

$$(2.11) \qquad |z - \beta_i| \sim |z - \alpha| + L_i(\alpha) \quad \text{for} \quad z \in D(\alpha) \ .$$

One inequality follows easily from $|z - \beta_i| \leq |z - \alpha| + |\beta_i - \alpha|$. For the reverse inequality, we have $|z - \beta_i| \geq |z - \alpha|$ (from the definition of $D(\alpha)$). Also, $L_i(\alpha) \sim |\beta_i - \alpha| \leq |z - \beta_i| + |z - \alpha| \leq 2|z - \beta_i|$. The estimate (2.11) is established.

For the remainder of the proof, we require the following two lemmas.

LEMMA 2.1.    Let $P(z) = \sum_{\nu=0}^{M} a_\nu z^\nu$, and let $\varepsilon, \lambda > 0$. Then

$$\int_I |P(e^{i\theta})|^\varepsilon d\theta \quad \sim \quad \sum_{\nu=0}^{M} |a_\nu|^\varepsilon$$

where $I \subseteq [0, 2\pi]$ is any interval of length at least $\lambda$, and the equivalence $\sim$ is up to constants depending only on $\varepsilon$, $\lambda$, but not on $P(z)$ and $I$ themselves.

*Proof of Lemma 2.1.* We may assume by homogeneity that $\sum |a_\nu|^\varepsilon = 1$. Then the integral is a continuous, nonvanishing function of the $a_\nu$ and of the two endpoints of $I$. Thus it is bounded above and below by positive constants.

LEMMA 2.2.     *Let $p, \delta \geq 0$, $N$ a positive integer and $c$ a real number such that $0 < c < 1$. Assume that $p - (N-k)\delta \neq 0$ for all integers $k$ such that $0 \leq k \leq N$. Then for every sequence of real numbers $\Lambda, L_0, \ldots, L_{N-1}$ satisfying $c\Lambda \geq L_0 \geq L_1 \cdots \geq L_{N-2} \geq L_{N-1} = 0$, the following estimate holds*:

$$(2.12) \quad \int_0^\Lambda \frac{r^p}{\prod_{i=0}^{N-1}(r+L_i)^\delta} \frac{dr}{r}$$

$$\sim \begin{cases} \left[ L_k^{(N-k)\delta-p} \prod_{0 \leq i < k} L_i^\delta \right]^{-1} & \text{if } (N-k-1)\delta < p < (N-k)\delta \\[2mm] \Lambda^{p-N\delta} & \text{if } p > N\delta. \end{cases}$$

*Here the equivalence $\sim$ is defined up to constants depending on $c$, $p$, $N$, and $\delta$, but not on $\Lambda$ and on the $L_i$, $0 \leq i \leq N-1$. In particular, each side of the estimate (2.12) is finite if and only if the other side is finite.*

*Proof of Lemma 2.2.* The lemma is evident if $L_0 = 0$, so we assume that $L_0 > 0$. Choose a constant $c_1$ with $1 < c_1 < c^{-1}$ and divide the interval of integration $(0, \Lambda)$ into $(0, c_1 L_0)$ and $(c_1 L_0, \Lambda)$. On the interval $(c_1 L_0, \Lambda)$, we have $r + L_i \sim r$ for all $i$, and we can write

$$(2.13) \qquad \int_{c_1 L_0}^\Lambda \frac{r^p}{\prod_{i=0}^{N-1}(r+L_i)^\delta} \frac{dr}{r} \sim \int_{c_1 L_0}^\Lambda r^{p-N\delta} \frac{dr}{r}.$$

If $p - N\delta > 0$, the right-hand side is of size $\Lambda^{p-N\delta}$. This establishes (2.12), since $\Lambda^{p-N\delta}$ is also an upper bound for the integral in (2.12) in this case. Thus we may assume that $p - N\delta < 0$, in which case (2.13) is of size $L_0^{p-N\delta}$. On the interval $(0, c_1 L_0)$, we can change scales $r \to L_0 r$ to obtain

$$(2.14) \qquad \int_0^{c_1 L_0} \frac{r^p}{\prod_{i=0}^{N-1}(r+L_i)^\delta} \frac{dr}{r} \sim L_0^{p-N\delta} \int_0^{c_1} \frac{r^p}{\prod_{i=0}^{N-1}(r+\frac{L_i}{L_0})^\delta} \frac{dr}{r}.$$

We can now argue by induction. We have already observed that the estimates hold in the case $N = 1$ (since in that case, we have $L_0 = 0$). Assume now that Lemma 2.2 holds for $N-1$. Since $c_1 \sim 1$, the right-hand side of (2.14) is of size

$$(2.15) \quad L_0^{p-N\delta} \int_0^{c_1} \frac{r^p}{\prod_{i=1}^{N-1}(r+\frac{L_i}{L_0})^\delta} \frac{dr}{r} = L_0^{p-N\delta} \int_0^{c_1} \frac{r^p}{\prod_{i=0}^{N-2}(r+\frac{L_{i+1}}{L_0})^\delta} \frac{dr}{r}.$$

The integrals in (2.15) are of the original form (2.12), with $N$ and $L_i$ replaced respectively by $N-1$ and $L_i^* = L_{i+1}/L_0$. The index $k$ in (2.12) gets replaced by $k-1$. Since $c_1 > 1 \geq L_i/L_0$, the induction hypothesis applies. We see immediately that (2.15) is of size $L_0^{p-N\delta}$ if $(N-1)\delta < p$, while for $p < (N-1)\delta$, it is of size

$$L_0^{p-N\delta}(L_{k-1}^*)^{p-(N-k)\delta} \prod_{0 \leq i < k-1} (L_i^*)^{-\delta} = L_0^{p-N\delta}(\frac{L_k}{L_0})^{p-(N-k)\delta} \prod_{1 \leq i < k} (\frac{L_i}{L_0})^{-\delta}$$
$$= L_k^{p-(N-k)\delta} \prod_{0 \leq i < k} L_i^{-\delta} \ ,$$

which is greater than $L_0^{p-N\delta}$. This proves Lemma 2.2.

It is now easy to establish the upper bounds in Theorem 1. By virtue of (2.11), the contribution from each region of integration $D(\alpha)$ can be estimated by

$$\int_{D(\alpha)} \frac{|P|^\varepsilon}{|Q|^\delta}dV \sim \int_{D(\alpha)} \frac{|P(z)|^\varepsilon}{\prod_{i=0}^{N-1}(|z-\alpha| + L_i(\alpha)|)^\delta} \ dV$$
$$< \int_0^{2\Lambda} \int_0^{2\pi} \frac{|P(\alpha + re^{i\theta})|^\varepsilon}{\prod_{i=0}^{N-1}(r+L_i(\alpha))^\delta} r d\theta dr,$$

where we have converted to polar coordinates centered at $\alpha$. Integrating with respect to $\theta$ and applying Lemma 2.1, we obtain

$$\int_{D(\alpha)} \frac{|P|^\varepsilon}{|Q|^\delta}dV \leq C \int_0^{2\Lambda} \frac{\sum_{\nu=0}^M |P^{(\nu)}(\alpha)|^\varepsilon r^{\nu\varepsilon+2}}{\prod_{i=0}^{N-1}(r+L_i(\alpha))^\delta} \ \frac{dr}{r}.$$

Applying Lemma 2.2 gives the upper bounds stated in Theorem 1.

*Lower bounds.* To establish the estimates in the other direction, fix $\alpha$ and, for $0 \leq \nu \leq N-2$, choose $r_\nu \geq 0$ and $\theta_\nu \in \mathbb{R}/\mathbb{Z}$ such that $(\beta_\nu - \alpha) = r_\nu e^{2\pi i \theta_\nu}$. Then $(\mathbb{R}/\mathbb{Z}) \setminus \{\theta_0, \dots, \theta_{N-2}\}$ is a disjoint union of intervals. Let $\psi$ be the midpoint of the largest interval (whose length is at least $1/N$). Let $\psi_0 = \psi - 1/4N$ and $\psi_1 = \psi + 1/4N$. Then

$$(2.16) \qquad \int_{B_\Lambda} \frac{|P(z)|^\varepsilon}{|Q(z)|^\delta} \ dV > \int_{\psi_0}^{\psi_1} \int_0^{\Lambda/2} \frac{|P(z)|^\varepsilon}{|Q(z)|^\delta} \ r dr d\theta$$

where $z = \alpha + re^{2\pi i\theta}$. Now for $z$ in the range $0 \leq r \leq 1/2$ and $|\theta - \psi| \leq 1/4N$, we have the estimate

$$(2.17) \qquad |z - \beta_\nu| \ \sim \ |\alpha - \beta_\nu| + |z - \alpha| \ \sim \ r + L_\nu(\alpha).$$

To see the first equivalence, we may, without loss of generality, assume that $\alpha = 0$. Then we simply observe that on the compact set $|z| + |\beta_\nu| = 1$, the function $|z - \beta_\nu|$ is continuous and positive, and thus bounded above and below by positive constants. The second equivalence follows from (2.10). Combining (2.16) and (2.17) and Lemma 2.1, we see that

$$\int_{B_\Lambda} \frac{|P(z)|^\varepsilon}{|Q(z)|^\delta} \, dV > \int_{\psi_0}^{\psi_1} \int_0^{\Lambda/2} \frac{|P(z)|^\varepsilon}{\prod_{i=0}^{N-1}(r + L_i(\alpha))^\delta} \, r dr d\theta$$

$$\sim \int_0^{\Lambda/2} \frac{\sum_{\nu=0}^M |P^{(\nu)}(\alpha)|^\varepsilon r^{\nu\varepsilon+2}}{\prod_{i=0}^{N-1}(r + L_i(\alpha))^\delta} \frac{dr}{r}.$$

Lemma 2.2 applies and gives the desired lower bounds for the integral of $|P(z)|^\varepsilon/|Q(z)|^\delta$. The proof of Theorem 1 is complete.

*Remarks.* (a) In the important special case where $P(z) = 1$ (and say, $N\delta - 2 > 0$), the estimate (2.6) reduces to

$$(2.18a) \qquad \int_{B_\Lambda} \frac{1}{|Q(z)|^\delta} \, dV \sim \sum_{\{\alpha; Q(\alpha)=0\}} \frac{1}{L_{k_0}(\alpha)^{(N-k_0)\delta-2} \prod_{0 \le i < k_0} L_i(\alpha)^\delta},$$

where the integer $k_0$ is defined by $(N - k_0 - 1)\delta < 2 < (N - k_0)\delta$. When $N\delta < 2$, the estimate (2.6) reduces to

$$(2.18b) \qquad \int_{B_\Lambda} \frac{1}{|Q(z)|^\delta} \, dV \sim \Lambda^{2-N\delta}.$$

(b) Another case of particular importance in this paper is the case when $Q(z)$ has no multiple roots; i.e., $L_k(\alpha) > 0$ for all $\alpha$ and all $k$, $0 \le k \le N - 2$. In this case, all expressions $\Phi_{\nu,k}(\alpha)$ as defined in (2.7) are nonvanishing for $0 \le l \le N - 2$. Furthermore, for each $\nu \ge 0$ and $k_\nu$ defined as in (2.5), it is readily verified that

$$(2.19) \qquad\qquad \Phi_{\nu,k}(\alpha) \ge \Phi_{\nu,k_\nu}(\alpha),$$

using the fact that the scales $L_k(\alpha)$ are decreasing in $k$. Thus the restricted sum over $(\nu, k_\nu)$ in Theorem 1 can be replaced unambiguously by the following sum over all $\nu$ and all $k$

$$(2.20) \qquad \int_{B_\Lambda} \frac{|P(z)|^\varepsilon}{|Q(z)|^\delta} \, dV \sim \sum_{\nu=0}^M \sum_{k=-1}^{N-2} \left\{ \sum_{\{\alpha; Q(\alpha)=0\}} \frac{|P^{(\nu)}(\alpha)|^\varepsilon}{\Phi_{\nu,k}(\alpha)} \right\}.$$

This expression has the advantage of being symmetric in the roots $\alpha$.

(c) Another suggestive form for the estimate (2.6), in the case where $Q(z)$ has no multiple roots, is the following

$$(2.21) \qquad \int_{B_\Lambda} \frac{|P(z)|^\varepsilon}{|Q(z)|^\delta} \, dV \sim \sum_{\{\alpha; Q(\alpha)=0\}} \sum_{k=-1}^{N-1} \frac{\sup_{|z-\alpha|=L_k(\alpha)} |P(z)|^\varepsilon}{L_k(\alpha)^{(N-k)\delta-2} \prod_{0 \le i < k} L_i(\alpha)^\delta},$$

valid e.g., when $\nu\varepsilon + 2 < N\delta$. (Other cases can also be expressed in the same way, with suitable modifications due to the expression for $\Phi_{\nu,k}(\alpha)$ when $\nu\varepsilon + 2 > N\delta$.) The expression (2.21) is a simple consequence of Theorem 1 and Lemma 2.1. It can be viewed as a natural generalization of (2.18) to the case of general $P(z)$.

(d) We had stated earlier that the assumption (2) in Theorem 1 can be removed by the inclusion of log terms. It is now evident that it suffices to incorporate such terms in Lemma 2.2 in the case where $\nu\varepsilon + 2 - (N - k)\delta$ vanishes. However, the resulting bounds would no longer belong to the class of "algebraic estimates."

## 3. Absolute cluster-scale estimates and symmetrization

Our next goal is to rewrite the cluster-scale estimates (2.6), and in particular the local cluster scales $L_k(\alpha)$ themselves, in terms of rational expressions in the coefficients of $P(z)$ and $Q(z)$. In the real setting of [18], this can only be done when $k < N/2$, in which case only the "absolute" cluster scales defined by

$$(3.1) \qquad L_k = \inf_\alpha L_k(\alpha)$$

mattered, and they can indeed be rewritten in terms of polynomials in the coefficients of $Q(z)$. In the present complex setting, it turns out that there are no such limitations, if we allow rational expressions in the coefficients of $Q(z)$. This is a first hint of significant differences between the two settings.

To be more specific, we need some additional notation: As before, we let $S$ be the set of zeroes of $Q$, counted with multiplicity. Note that the absolute cluster scales $L_k$ defined by (3.1) are also given by

$$L_k = \min_{S_{N-k}} \{d(S_{N-k})\}$$

where the minimum is taken over all subsets $S_{N-k} \subseteq S$ with $|S_{N-k}| = N - k$. Let $a_i$ be the coefficients of $Q(z)$:

$$(3.2) \qquad Q(z) = \prod_{\alpha \in S}(z - \alpha) = \sum_{i=0}^{N} a_i y^{N-i},$$

and introduce for each integer $r$ with $1 \le r \le N/2$, the $r$-discriminant:

$$(3.3) \qquad \Delta_r(\alpha_1, \alpha_2, \ldots, \alpha_N) = \Delta_r = \sup_M \prod_{\nu=1}^{r} |\alpha_{i_\nu} - \alpha_{j_\nu}|,$$

where the supremum is taken over all $2r$-tuples $M = (i_1, \ldots, i_r, j_1, \ldots j_r)$ consisting of distinct integers with $1 \le i_\nu, j_\nu \le N$ for all $\nu$. Such an $2r$-tuple is said to be *admissible*. We say that an admissible $2r$-tuple of distinct integers is maximizing for the set $S$ if it achieves the supremum in (3.3).

The following theorem relates cluster scales of a monic polynomial $Q(z)$ to polynomial expressions in its coefficients:

THEOREM 2.    (i) *There is a root $\alpha \in S$ with*

$$(3.4) \qquad L_i(\alpha) \sim L_i, \quad \text{for all integers } i, \ 0 \leq i \leq N/2.$$

(ii) *For all integers $r$ satisfying $1 \leq r \leq N/2$,*

$$(3.5) \qquad \Delta_r(\alpha_1, \ldots, \alpha_N) \sim L_0 L_1 \cdots L_{r-1}.$$

(iii) *Let $h = N!/(2r)!$. There are polynomials*

$$D_{r,1}, D_{r,2}, \ldots D_{r,h} \in \mathbb{Z}[A_1, \ldots, A_N]$$

*such that*

$$(3.6) \qquad \Delta_r(\alpha_1, \ldots, \alpha_N) \sim \left\{ \sum_{q=1}^{h} |D_{r,q}(a_1, \ldots, a_N)| \right\}^{1/h!}.$$

(iv) *For each $k$, $1 \leq k \leq N-2$, there exist polynomials*

$$\sigma_{k,i} \in \mathbb{Z}[A_1, \ldots, A_N; Z]$$

*with polynomial coefficients in $a_1, \cdots, a_N$, so that*

$$(3.7) \qquad \prod_{0 \leq i \leq k} L_i(\alpha) \sim \sum_{i=1}^{N!/(N-k-1)!} |\sigma_{k,i}(a_1, \cdots, a_N; \alpha)|^{\frac{1}{i}}.$$

(Here all equivalences are defined up to constants depending only on $N$.)

*Proof of Theorem 2.* We begin with the proof of (i). For each $i \leq k = [N/2]$, choose $S_{N-i} \subset S$ such that $|S_{N-i}| = N-i$ and $d(S_{N-i}) = L_i$. We claim that any $\alpha \in S_{N-k}$ satisfies the required property. Indeed, for any $i \leq k$, $S_{N-i}$ and $S_{N-k}$ have nonempty intersection. This is evident if $i = k$, and if $i < k$, it follows from the fact that $(N-i)+(N-k) > N$ (note that $i+k < 2k \leq N$, which implies that $i < N-k$). Hence we must have $d(S_{N-i} \cup S_{N-k}) \sim d(S_{N-i}) = L_i$. On the other hand, the set $S_{N-i} \cup S_{N-k}$ has at least $N-i$ elements and it contains $\alpha$. Thus $d(S_{N-i} \cup S_{N-k}) \geq L_i(\alpha) \geq L_i$. This shows that $L_i(\alpha) \sim L_i$, and (i) is proved.

Next, we observe that there exists a partition of $S$ into two disjoint proper subsets $S_{N_1}, S_{N_2}$, $N_1, N_2 \geq 1$, with $|S_{N_i}| = N_i < N$ and $N_1 + N_2 = N$, and $|\alpha_1 - \alpha_2| \geq L_0/N$ for all $\alpha_1, \alpha_2$ with $\alpha_i \in S_{N_i}$. To see this, we define an equivalence relation on the set $S$ as follows: $\alpha \sim \alpha'$ if there exist $\beta_1, \ldots, \beta_m \in S$ such that $\alpha = \beta_1$, $\alpha' = \beta_m$ and $|\beta_i - \beta_{i+1}| < L_0/N$ for all $i$. Fix $\alpha \in S$ and let $S_{N_1}$ be the equivalence class containing $\alpha$, and $S_{N_2}$ the complement of $S_{N_1}$. This is the desired partition.

We may assume that $N_1 \leq N_2$. Then we have the following equivalences:

$$(3.8) \qquad L_0 \sim L_1 \sim \cdots \sim L_{N_1-1}.$$

Since $2r \leq N$ we must have $r \leq N_2$. We consider two cases:

*Case* 1. $r \leq N_1 \leq N_2$. It is clear that $\Delta_r \leq L_0^r$. Since $r \leq N_1$, we can choose an admissible $2r$-tuple $M$ such that $\alpha_{i_1}, \alpha_{i_2}, \ldots, \alpha_{i_r} \in S_{N_1}$ and $\alpha_{j_1}, \alpha_{j_2}, \ldots, \alpha_{j_r} \in S_{N_2}$. This shows that $\Delta_r \sim L_0^r$. Thus (3.8) implies (ii) in this case.

*Case* 2. $N_1 < r \leq N_2$. Let $M$ be a maximizing admissible $2r$-tuple for $S$. We may assume that if $\{i_\nu, j_\nu\} \cap S_{N_1}$ is nonempty for some $\nu$, then $i_\nu \in S_{N_1}$ (interchange $i_\nu$ and $j_\nu$ if necessary). We may further assume that $S_{N_1} \subseteq \{i_1, i_2, \ldots, i_r\}$: If not, then there exists $i \in S_{N_1}$ such that $i \neq i_\nu$ for all $\nu$. Then $i \neq j_\nu$ for all $\nu$ (by assumption). Choose $\nu$ such that $i_\nu \in S_{N_2}$ (such a $\nu$ exists since $r > N_1$). Then $j_\nu \in S_{N_2}$ as well (by assumption). Since $|\alpha_i - \alpha_{j_\nu}| \sim L_0$, we can replace $i_\nu$ by $i$ without decreasing the size of $\Delta_r$, while increasing $|S_{N_1} \cap \{i_1, \ldots, i_r\}|$. Continuing in this fashion, we see that we may assume that $S_{N_1} \subseteq \{i_1, i_2, \ldots, i_r\}$ which implies that

$$\Delta_r \leq C \cdot L_0^{N_1} \cdot \Delta_{r-N_1}(S_{N_2}).$$

To get the inequality in the reverse direction, observe that $N_2 - 2(r - N_1) \geq N_1$ so that when we choose a maximizing admissible $2(r - N_1)$-tuple for the set $S_{N_2}$, there are at least $N_1$ elements left over which can be paired with the $N_1$ elements in the set $S_{N_1}$. Thus we have proved

$$(3.9) \qquad \Delta_r \sim L_0^{N_1} \cdot \Delta_{r-N_1}(S_{N_2}).$$

Since $2(r - N_1) \leq N_2$, we can use induction to deduce:

$$(3.10)$$
$$\Delta_{r-N_1}(S_{N_2}) \sim L_0(S_{N_2})L_1(S_{N_2})\cdots L_{(r-N_1-1)}(S_{N_2}) \sim L_{N_1}L_{N_1+1}\cdots L_{r-1}.$$

Combining (3.8), (3.9) and (3.10) we obtain (ii).

To prove (iii), we make use of the following elementary fact:

$$(3.11) \qquad \sum_{i=1}^{n} |\gamma_i| \ \sim \ \sum_{r=1}^{n} \Big| \sum_{i=1}^{n} \gamma_i^r \Big|^{1/r}$$

for any $n$ complex numbers $\gamma_1, \cdots, \gamma_n$. Define then the polynomial $F_r(T)$ by

$$(3.12) \quad F_r(T) = \prod_M \Big( T - \prod_{\nu=1}^{r}(\alpha_{i_\nu} - \alpha_{j_\nu}) \Big) \ = \ \sum_{q=0}^{H} B_{r,q}(a_1, \ldots, a_N) T^{H-q}$$

where the product is taken over all admissible $2r$-tuples $M$. The coefficients $B_{r,q}$ of $F_r(T)$ are in $\mathbb{Z}[A_1, \ldots, A_N]$. Now $\Delta_r$ is the size of the largest root of $F_r$, and hence is of the size given by the right-hand side of (3.11), with $\gamma_i$ the roots of $F_r(T)$. But symmetric polynomials in the $\gamma_i$'s are polynomials in the coefficients $B_{r,q}(a_1, \cdots, a_N)$, and hence polynomials in the $(a_1, \cdots, a_N)$ themselves. Clearly, they are of the form described in (3.6), and (iii) is proved.

We turn now to the proof of (iv). Let $0 \leq k \leq N - 2$. Then we note that

$$(3.13) \quad L_0(\alpha) L_1(\alpha) \cdots L_k(\alpha) \ \sim \ \sup_{\lambda \in \Lambda_k} \prod_{\nu=0}^{k} |\alpha - \alpha_{i_\nu}| \sim \sum_{\lambda \in \Lambda_k} \prod_{r=0}^{k} |\alpha - \alpha_{i_\nu}|$$

where $\Lambda_k = \{\lambda = (\alpha_{i_0}, \ldots, \alpha_{i_k}) : 1 \leq i_\nu \leq N, \ \ i_\nu \neq i_\mu \ \text{if} \ \ \mu \neq \nu\}$. To see (3.13), we select successively $\beta_0, \cdots, \beta_{N-1}$ as in the proof of Theorem 1 so that $L_i(\alpha) \sim |\beta_i - \alpha|$ for $0 \leq i \leq N-1$ (c.f. (2.10)). This implies that the left-hand side of (3.13) is bounded by the right-hand side. To see the opposite inequality, we note that if $(\alpha_{i_0}, \cdots, \alpha_{i_k})$ is any sequence appearing in the sup on the right-hand side, and if we order them in decreasing order of their distances to $\alpha$, say $|\alpha_{i_0} - \alpha| \geq |\alpha_{i_1} - \alpha| \geq \cdots \geq |\alpha_{i_k} - \alpha|$, then $L_j(\alpha) \geq |\alpha_{i_j} - \alpha|$ for $j \leq k$. The estimate (3.13) follows.

For $\lambda = (\alpha_{i_0}, \cdots, \alpha_{i_k}) \in \Lambda_k$, let $G_\lambda(T) = \prod_{\nu=0}^{k}(T - \alpha_{i_\nu})$. Let $\sigma_{k,1}(T)$, $\sigma_{k,2}(T), \ldots \sigma_{k,|\Lambda_k|}(T)$ be the standard symmetric polynomials in the $G_\lambda$; i.e.

$$(3.14) \quad \sigma_{k,i} = \sum_{\lambda \in \Lambda_k} G_\lambda^i, \ \ 1 \leq i \leq |\Lambda_k|.$$

Then $\sigma_{k,i} = \sigma_{k,i}(a_1, a_2, \ldots, a_N; T)$ where

$$\sigma_{k,i}(A_1, A_2, \ldots, A_N; T) \in \mathbb{Z}[A_1, \ldots, A_N, T].$$

Thus (3.11) and (3.13) imply that

$$(3.15) \quad L_0(\alpha) L_1(\alpha) \cdots L_k(\alpha) \ \sim \ \sum_{i=1}^{|\Lambda_k|} |\sigma_{k,i}(a_1, a_2, \ldots, a_N; \alpha)|^{1/i}$$

and the proof of Theorem 2 is complete.

As a first step in the program of rewriting all integrals of the form (2.6) in terms of the coefficients of $P(z)$ and $Q(z)$, we focus in the remaining part of this section on the special case where the numerator $P(z)$ is the constant 1. As in Theorem 2 where the case $k \leq N/2$ is significantly simpler than the case of general $k$, here the corresponding case $\delta < 4/N$ is significantly simpler than the case of general $\delta$. This fact is at the origin of the considerably greater difficulties which arise in the treatment of stability in dimensions $n \geq 3$, compared to dimensions $n \leq 2$. We refer to Sections 5–7, and especially Sections 5.B and 6 for a fuller discussion.

THEOREM 3. (i) *Let $N$ be nonnegative integer, and let $\delta$ be a positive rational number with $\frac{2}{\delta} \notin \{1, \cdots, N\}$. Then there exist integers $N'$, $I$, and $J$, rational numbers $\delta'$ and $\varepsilon'$, and polynomials $F_1, \cdots, F_I$, $G_1, \cdots, G_J \in Z[A_1, \cdots, A_N]$ with degrees bounded by $N'$, so that*

$$(3.16) \qquad \int_{B_\Lambda} |Q(z)|^{-\delta} \; dV \; \sim \; \frac{(\sum_{i=1}^I |F_i(a_1, \cdots, a_N)|^2)^{\varepsilon'/2}}{(\sum_{j=1}^J |G_j(a_1, \cdots, a_N)|^2)^{\delta'/2}},$$

*for all monic polynomials $Q(z) = \sum_{i=0}^N a_i z^{N-i}$, whose zeroes all lie inside the ball $B_{\Lambda/2}$. The equivalence in (3.16) means the following. The left-hand side is infinite if and only if $\sum_{j=1}^J |G_j(a_1, \cdots, a_N)| = 0$. When this is not the case, both sides are finite and bounded by each other, up to constants depending only on $N$, $\delta$, and $\Lambda$.*

(ii) *When $\delta < 4/N$, we can take $I = 1$ and $F_1(z) = 1$ in the expression (3.16).*

*Proof of Theorem* 3. We apply Theorem 1. When $N\delta - 2 < 0$, Theorem 1 implies that the integral under consideration is of size $\Lambda^{2-N\delta}$, which is obviously of the desired form (3.16). Otherwise, we use the form (2.18) of Theorem 1. We can write

$$(3.17) \qquad \begin{aligned} \Phi_{0,k_0}(\alpha) &= L_{k_0}(\alpha)^{(N-k_0)\delta-2} \prod_{0 \le i < k_0} L_i(\alpha)^{-\delta} \\ &= \big[ \prod_{0 \le i \le k_0 - 1} L_i(\alpha) \big]^{\varepsilon_1} \big[ \prod_{0 \le i \le k_0} L_i(\alpha) \big]^{\varepsilon_2}, \end{aligned}$$

where $\varepsilon_2 = (N - k_0)\delta - 2$ and $\varepsilon_1 = 2 - (N - k_0 - 1)\delta$ are positive numbers. Since each factor on the above right-hand side can be expressed in the form (iv) of Theorem 2, and since both $\varepsilon_1$ and $\varepsilon_2$ are rational numbers, it follows that the size of $\Phi_{0,k_0}(\alpha)$ can be expressed in the form

$$(3.18) \qquad \Phi_{0,k_0}(\alpha) \sim \big[ \sum_{j=1}^\eta |K_j(a_1, \cdots, a_N; \alpha)|^2 \big]^{\frac{\varepsilon_3}{2}},$$

where $K_j(a_1, \cdots, a_N; \alpha)$ are polynomials in all variables and $\varepsilon_3$ is a rational number. The integral in (3.16) is infinite if and only if $\Phi_{0,k_0}(\alpha) = 0$ for some $\alpha$. The sum

$$(3.19) \qquad \begin{aligned} &\sum_{\{\alpha; Q(\alpha)=0\}} \frac{1}{\big[ \sum_{j=1}^\eta |K_j(a_1, \cdots, a_N; \alpha)|^2 \big]^{\frac{\varepsilon_3}{2}}} \\ &\qquad \sim \Big\{ \sum_{\{\alpha; Q(\alpha)=0\}} \frac{1}{\sum_{j=1}^\eta |K_j(a_1, \cdots, a_N; \alpha)|^2} \Big\}^{\frac{\varepsilon_3}{2}} \end{aligned}$$

can be reduced to the same denominator

$$\prod_{\{\alpha;\, Q(\alpha)=0\}} \sum_{j=1}^{\eta} |K_j(a_1,\cdots,a_N;\alpha)|^2.$$

This denominator vanishes if and only if one of the factors vanishes, i.e., if and only if the integral in (3.16) is infinite. Furthermore, multiplying it out gives a sum of $n$ terms, $\gamma_1,\cdots,\gamma_n$. Applying the elementary symmetry principle, it can be rewritten in the form of the denominator appearing in (3.16) of Theorem 3. The numerator of (3.16) is obtained by similar arguments. The result is clearly of the form (3.16) and (i) is established.

To establish (ii), we note that in the range $\delta < 4/N$, we have $k_0 \leq N/2$. Thus (i) of Theorem 2 shows that there exists a root $\alpha$ with $L_i(\alpha) \sim \min_\beta L_i(\beta) = L_i$, for all $i \leq k_0$. This implies that

$$(3.20)$$

$$\sum_{\{\alpha;\, Q(\alpha)=0\}} \Phi_{0,k_0}^{-1}(\alpha) \sim L_k^{2-(N-k)\delta} \prod_{0 \leq i < k_0} L_i^{-\delta} = \Big[\prod_{0 \leq i \leq k_0-1} L_i\Big]^{-\varepsilon_1} \Big[\prod_{0 \leq i \leq k_0} L_i\Big]^{-\varepsilon_2}.$$

In view of (iii) of Theorem 2, the right-hand side of (3.20) is of the desired form.

## 4. Algebraic estimates

A. *Stratification of spaces of absolute rational powers.* In Theorem 1, we derived a sharp formula for integrals of rational functions $|P(z)|^\varepsilon |Q(z)|^{-\delta}$ in terms of local scales $L_k(\alpha)$ for $Q(z)$. In Theorem 3, we have seen how, in some relatively simpler cases, these formulas can be rewritten in terms of ratios of absolute values of polynomial expressions in the coefficients of $P(z)$ and $Q(z)$. It is now important for us to extend these results in several directions:

(1) In order to treat functions of several variables, we need to iterate the estimates, and hence extend our analysis from integrals of rational expressions $|P(z)|^\varepsilon |Q(z)|^{-\delta}$ to integrals of the form

$$(4.1) \qquad \int_{B_\Lambda} \frac{(\sum_{i=1}^I |P_i(z)|^2)^{\varepsilon/2}}{(\sum_{j=1}^J |Q_j(z)|^2)^{\delta/2}}\, dV$$

where $P_i(z)$ and $Q_j(z)$ are polynomials. We shall refer to the expression which appears as the integrand in (4.1) as an ARP, which is an acronym for *absolute rational power.* Compared with the case of real-analytic functions discussed in [18], the ARP's are a source of significant additional complications.

(2) We must also rewrite, in complete generality, the cluster scales estimates (2.6) in terms of the coefficients of $P_i(z)$ and $Q_j(z)$. This "symmetrization" was possible when $I = J = 1$ and $P_1(z) = 1$ (cf. Theorem 3). We can symmetrize as well in the case where $J = 1$, and the polynomial $Q_1(z)$ has only simple roots. In this latter case, formula (2.20) gives an expression for the size of the integral which is manifestly symmetric in the roots $\alpha$ of $Q(z)$. Such formulas can be rewritten in terms of the coefficients of $Q(z)$, as we will show in greater detail below. The difficulty in the general case is that the polynomials $Q_j(z)$ in the denominator may vanish of high and varying orders at each of their roots, with uneven compensating effects from the behavior of the polynomials $P_i(z)$ near these roots. Clearly, such phenomena can lead to high instability, and must be treated with some care.

It turns out that there is a simple and cogent geometric description for the structure of the general symmetrization process. The precise formulation will be given in Theorem 4, but the basic idea is the following: Once we fix $\varepsilon, \delta$ and an upper bound on the degrees of the polynomials, the space $\mathbb{C}|(Z)|$ of ARP's, i.e., of functions of the form $((\sum_{i=1}^{I} |P_i(z)|^2)^{\varepsilon/2}/((\sum_{j=1}^{J} |Q_j(z)|^2)^{\delta/2})$, may be regarded, via the coefficients of the $P_i$ and $Q_j$, as an affine space over the complex numbers of finite dimension. We shall show that this ARP space admits a *stratification by a finite number of algebraic varieties*. This stratification induces, in a canonical way, a decompositon of the ARP space into a finite disjoint union of quasi-affine pieces. On each quasi-affine piece the size of the integral of the ARP is controlled by a corresponding expression in $\mathbb{Z}|(b)|$, i.e., an ARP in in the coefficients $b$ of $P_i(Z)$ and $Q_j(Z)$, with integer coefficients.

To state our results more precisely, we formulate more carefully the above concept of ARP, and introduce some notation.

*Definition* 4.1. Let $N$ be a positive integer and let $A$ be a subring of $\mathbb{C}$. Then an absolute rational power (an "ARP") in $N$ variables with coefficients in $A$ is an expression of the form

$$(4.2) \quad R(Z) = R(Z_1, \ldots, Z_N) = \frac{(\sum_{i=1}^{I} |P_i(Z_1, \ldots, Z_N)|^2)^{\varepsilon/2}}{(\sum_{j=1}^{J} |Q_j(Z_1, \ldots, Z_N)|^2)^{\delta/2}} = \frac{K(Z)}{L(Z)}.$$

Here $Z = (Z_1, \ldots, Z_N)$ is a vector of independent variables, $i$ and $j$ range over finite index sets, the $P_i$ and the $Q_j$ are polynomials in $N$ variables with coefficients in $A$, $\varepsilon$ and $\delta$ are nonnegative rational numbers, and at least one of the $P_i$ is nonzero. The numerator $K(Z)$, and denominator $L(Z)$ will be called absolute polynomial powers.

We shall sometimes write $R = R_b$, where $b = (b', b'')$ and $b'$, $b''$ are respectively the vector of coefficients of all the $P_i(Z)$'s and the $Q_j(Z)$'s. Thus $b$ lies is a certain finite-dimensional free $A$ module whose dimension is bounded in terms of the maximum of the degrees of $P_i(Z)$ and $Q_j(Z)$, and the number $I + J$ of terms which appear in the numerator and denominator of (4.2).

The set of all ARP's in $N$ variables will be denoted by $A|(Z)| = A|(Z_i, \ldots, Z_N)|$. It is clearly closed, up to size, under addition, multiplication, and division by nonzero elements.

The spaces of ARP's of particular interest in this paper are the following:

- $\mathbb{Z}|(Z_1, \cdots, Z_N)| = \mathbb{Z}[Z_1, \cdots, Z_{N-1}]|(Z_N)|$: ARP's in $N$ variables with integer coefficients, which can also be thought of as ARP's in the variable $Z_N$, with coefficients in the ring of polynomials in $(Z_1, \cdots, Z_{N-1})$ with integer coefficients;

- $\mathbb{C}|(Z_1, \cdots, Z_N)| = \mathbb{C}[Z_1, \cdots, Z_{N-1}]|Z_N|$: ARP's in $N$ variables with complex coefficients;

- $H\{Z_1, \cdots, Z_{N-1}\}|(Z_N)|$: ARP's in the variable $Z_N$, with coefficients which are holomorphic functions in $(Z_1, \cdots, Z_{N-1})$.

Each $R(Z_1, \ldots, Z_N) = K(Z_1, \ldots, Z_N)/L(Z_1, \ldots, Z_N) \in \mathbb{C}|Z_1, \cdots, Z_N|$ corresponds to a continuous $(0, \infty]$ valued function $R(z_1, \ldots, z_N)$ which is defined on the "domain of $R$," that is the complement of the affine subvariety of $\mathbb{C}^N$ defined by $K(z) = 0$. The "strict domain of $R$" is defined to be the complement of the affine variety defined by $K(z) = L(z) = 0$, so that if we restrict $R$ to its strict domain, it takes values in $(0, \infty)$. As in the case of polynomials, we shall often use upper cases letters vs. lower cases letters to distinguish between the formal expression $R(Z_1, \cdots, Z_N)$ and its realization $R(z_1, \cdots, z_N)$ as a function on its domain of definition.

The following simple observations are also important in the sequel:

1. Since both $\varepsilon$ and $\delta$ are rational, the ARP $R(z)$ is always equivalent in size to another ARP $R^{\#}(z)$ of the form

$$(4.2a) \qquad R^{\#}(z) = \left( \frac{\sum_{i=1}^{I} |P_i^{\#}(Z_1, \ldots, Z_N)|^2}{\sum_{j=1}^{J} |Q_j^{\#}(Z_1, \ldots, Z_N)|^2} \right)^{\varepsilon^{\#}/2}.$$

Indeed, we can just write $\varepsilon = A/D$, $\delta = B/D$ with $A, B, D \in \mathbb{N}$, and set $P_i^{\#}(Z) = (P_i(Z))^A$, $Q_j^{\#}(Z) = (Q_j(Z))^B$ and $\varepsilon^{\#} = 1/D$. When only the size of $R(z)$ matters, we shall freely replace $R(z)$ by $R^{\#}(z) \sim R(z)$. Similarly, we shall also often replace $R(z)$ by $\sum_{i=1}^{I} |P_i(z_1, \cdots, z_N)|^{\varepsilon} / \sum_{j=1}^{J} |Q_j(z_1, \cdots, z_N)|^{\delta}$ without additional comment.

2. If $R(Z) \in \mathbb{C}|(Z)|$ is a one-variable ARP, then its domain of definition is the complement of a finite set of points. We claim that $R(z)$ has a unique extension to a continuous function on all of $\mathbb{C}$: If $z_0$ is a point outside the domain of definition of $R$, then $K(z_0) = L(z_0) = 0$. This means that $|z - z_0|$ is a factor of both the numerator and denominator. Cancelling out the highest power of $|z - z_0|$ which divides both $K(z)$ and $L(z)$, we see that we can extend $R(z)$ to a continuous function in a neighborhood of $z_0$.

We can now state the main theorem on algebraic estimates, which allows us to calculate the size of a one dimensional ARP integral: Fix nonnegative rational numbers $\varepsilon$, $\delta$ and let $I, J, M_i, N_j$, be integers with $N_1 \geq N_2 \geq \cdots \geq N_J$. Let $\delta$ be a nonnegative rational number such that $(\varepsilon, \delta)$ is nondegenerate; that is, assume that

$$(4.3) \qquad \nu\varepsilon + 2 - \frac{l}{[\delta] + 1}\delta \neq 0,$$

for any integer $0 \leq \nu \leq \max M_i$ and any integer $0 \leq l \leq ([\delta] + 1)\max N_i$. Let $C > 0$ be a real number. Define

$$(4.4) \quad \mathcal{B} = \{b = (b'_{i,\mu}, b''_{j,\nu}) : 1 \leq i \leq I, 1 \leq j \leq J, 0 \leq \mu \leq M_i,$$
$$0 \leq \nu \leq N_j, b'_{i,\mu}, b''_{j,\nu} \in \mathbb{C}, b'_{i,\mu} \neq 0 \text{ for some } i, \mu\}.$$

Thus $\mathcal{B} \subseteq \mathbb{C}^a$ for some $a$. To every point $b = (b', b'') \in \mathcal{B} = \mathcal{B}' \times \mathcal{B}''$ we can associate a one-variable ARP $R_b \in \mathbb{C}|(Z)|$ as follows:

$$(4.5) \qquad R_b(Z) = \frac{(\sum_{i=1}^I |\sum_{\mu=0}^{M_i} b'_{i,\mu} Z^\mu|^2)^{\varepsilon/2}}{(\sum_{j=1}^J |\sum_{\nu=0}^{N_j} b''_{j,\nu} Z^\nu|^2)^{\delta/2}} = \frac{(\sum_{i=1}^I |P_i(b'; Z)|^2)^{\varepsilon/2}}{(\sum_{j=1}^J |Q_j(b''; Z)|^2)^{\delta/2}}.$$

We shall also use the following norm for polynomials of a bounded degree. For $Q(Z) = \sum_{\mu=0}^M b_\mu Z^\mu \in \mathbb{C}[Z]$, let

$$(4.6) \qquad |||Q||| = \sum_{\mu=0}^M |b_\mu|.$$

Finally, we recall that an algebraic variety in $\mathbb{C}^a$ is a subset which can be realized as the simultaneous vanishing of a finite set of polynomials in $\mathbb{C}[X_1, \ldots, X_a]$. In particular, we do not assume that our varieties are irreducible.

THEOREM 4.    *There exist a real number* $s \in (0, 1)$, *a finite chain* $\mathcal{U}_\lambda$, $0 \leq \lambda \leq \mathcal{N}$, *of algebraic varieties*

$$(4.7) \qquad \mathcal{B} = U_0 \supset U_1 \supset \cdots \supset U_\mathcal{N} = \emptyset \ ,$$

*and a corresponding sequence* $T_\lambda \in \mathbb{Z}[B'_{i,\mu}, B''_{j,\nu}]$, $0 \leq \lambda \leq \mathcal{N} - 1$, *with the following properties*:

1.  *If* $T_\lambda(B) = K_\lambda(B)/L_\lambda(B)$ *as in* (4.2), *and if* $b \in U_\lambda \backslash U_{\lambda+1}$, *then* $K_\lambda(b) \neq 0$.
    *In particular,* $T_\lambda$ *is defined, and nowhere vanishing on* $U_\lambda \backslash U_{\lambda+1}$.

2. *Fix $C > 0$. For $b \in U_\lambda \backslash U_{\lambda+1}$ satisfying $|||Q_1(b; Z) - Z^{N_1}||| < s$ and $|||Q_j(b; Z)||| < C$ for $j > 1$, the vector $b$ is in the domain of $T_\lambda$ and*

$$(4.8) \qquad\qquad \int_B R_b(z) \ dV \sim T_\lambda(b).$$

*Here $B$ is the unit disk in $\mathbb{C}$. The implied constants depend only on the fixed constants $I, J, M_i, N_j, \varepsilon, \delta, C$, and are independent of $b$;*

3. *We have*

$$(4.9)$$
$$\mathcal{U}_{\lambda+1} = \{ b \in U_\lambda \ : \ K_\lambda(b) = 0 \} \ = \ \{ b \in U_\lambda \ : \ K_\lambda(b) = 0 \text{ and } L_\lambda(b) = 0 \}.$$

*Remark 1.* Let $A = [\delta] + 1$. Then replacing $Q_j$ by $Q_j^A$ and $\delta$ by $\delta/A$, we see that it suffices to prove the theorem in the case $\delta < 1$.

*Remark 2.* The theorem implies that $\int_B R_b(z) \ dV = \infty$ if and only if $b$ lies in a certain constructible set, that is, a set which can be described by taking a finite number of intersections, unions and complements of algebraic varieties.

We shall make use of the following notation: If $b \in \mathcal{B}$, then we shall define

$$(4.10) \qquad\qquad\qquad \lambda(b) = \lambda$$

to be the unique integer $\lambda$ for which $b \in U_\lambda \backslash U_{\lambda+1}$.

The proof of Theorem 4 will be the subject of Sections B through G. It requires several types of arguments, which we present separately for the sake of clarity.

B. *Symmetrization in the roots and absolute rational powers.* We have seen from Theorem 1 that the integration of rational functions $|P(z)|^\varepsilon |Q(z)|^{-\delta}$ gives rise to sums, over the roots, $\alpha$, of $Q(z)$, of the expressions $|P^\nu(\alpha)|^\varepsilon \Phi_{\nu, k_\nu}(\alpha)^{-1}$. The term $\Phi_{\nu, k_\nu}(\alpha)$ is originally defined in terms of the local scales $L_k(\alpha)$ of $Q(z)$. However, in view of part (iv) of Theorem 2, these expressions can all be rewritten as ARP's in $\mathbb{Z}[b', b'']|\alpha|$, i.e., as ARP's in $\alpha$, with coefficients which are polynomials in the coefficients of $P(z)$ and $Q(z)$. We need to address now the problem of rewriting the sum over $\alpha$ of such expressions as the value, at $(b', b'')$, of an ARP in $\mathbb{Z}|B', B''|$. In the simplest case, this can be done algebraically through the following lemmas.

LEMMA 4.1. *Fix positive integers $M_1, \ldots, M_I, N_1, \ldots, N_J$ and $d$. Fix rational numbers $\varepsilon, \delta \geq 0$. Then there exists $T \in \mathbb{Z}|(A, B)|$ where $A = (A_0, \ldots, A_d)$, $B = (B'_{i,\mu}, B''_{j,\nu})$, $1 \leq i \leq I$, $1 \leq j \leq J$, $0 \leq \mu \leq M_i$*

and $0 \leq \nu \leq N_j$ with the following properties: For every one-variable ARP $R = R_b \in \mathbb{C}|Z|$ of the form

$$(4.11) \qquad R_b(Z) = \frac{\sum_{i=1}^{I} |\sum_{\mu=0}^{M_i} b'_{i,\mu} Z^\mu|^\varepsilon}{\sum_{j=1}^{J} |\sum_{\nu=0}^{N_j} b''_{j,\nu} Z^\nu|^\delta} = \frac{\sum_{i=1}^{I} |P_i(Z)|^\varepsilon}{\sum_{j=1}^{J} |Q_j(Z)|^\delta}$$

and for every polynomial $G(u) = a_0 u^d + a_1 u^{d-1} + \cdots + a_N \in \mathbb{C}[u]$ of degree $d$, whose roots all lie in the strict domain of $R_b(z)$, we have

1. The point $(a, b)$ lies in the strict domain of $T$, where $a = (a_0, \ldots, a_n)$ and $b = (b'_{i,\mu}, b''_{j,\nu})$. Moreover, $T(a, b) \in (0, \infty)$.

2. The following estimate holds:

$$(4.12) \qquad \sum_{\{\alpha : G(\alpha) = 0\}} R_b(\alpha) \sim T(a, b)$$

and the implied constants depend only on $M_i, N_j, d, \varepsilon, \delta$.

Note that implicit in the assumptions of the lemma is the condition $a_0 \neq 0$.

Lemma 4.1 is a simple consequence of the following two lemmas, which generalize the arguments used earlier in the proof of Theorem 3.

Let L be a commutative ring, $n$ a positive integer, and let $F \in$ L$[X_1, \ldots, X_n]$ be a polynomial in $X = (X_1, \ldots, X_n)$. If $\pi \in S_n$ is a permutation of the integers $\{1, 2, \ldots, n\}$, we write $X^\pi = (X_{\pi(1)}, \ldots, X_{\pi(n)})$. We say that $F$ is symmetric if $F(X^\pi) = F(X)$ for all $\pi \in S_n$.

For $1 \leq j \leq n$ we define the standard symmetric polynomials, $\sigma_j \in \mathbb{Z}[X_1, \ldots, X_n]$, by the equation

$$\prod_{i=1}^{n} (T + X_i) = T^n + \sum_{j=1}^{n} \sigma_j T^{n-j}.$$

The following lemma is well-known (for example, a proof may be found in [13]):

Lemma 4.2. Let $F \in L[X]$ be a polynomial in $n$ variables. Then $F$ is symmetric if and only if $F \in L[\sigma_1, \ldots, \sigma_n]$, that is, if and only if $F(X) = \tilde{F}(\Sigma)$ where $\Sigma = (\sigma_1, \ldots, \sigma_n)$ and $\tilde{F}$ is a polynomial in $n$ variables with coefficients in $L$.

The next lemma generalizes Lemma 4.2 to the case where $F$ is an absolute polynomial power (that is, an ARP whose denominator is one).

Lemma 4.3. Fix $D > 0$. Let $F \in C|(A; X)|$ be an absolute polynomial power in the variables $A = (A_1, \ldots, A_m)$ and $X = (X_1, \ldots, X_n)$, where $C$ is a subring of $\mathbb{C}$. Let all polynomials occurring in $F$ have degrees bounded by $D$. Assume that $F$ is symmetric in the $X$ variables; that is, *assume $F(A; X^\pi) = F(A; X)$ for all $\pi \in S_n$. Then $F(A; X) \sim \tilde{F}(A; \Sigma)$ for some absolute polynomial power $\tilde{F} \in C|(A; \Sigma)|$. The implied constants depend only on $\varepsilon, \delta$ and $D$.*

*Proof of Lemma* 4.3.

*Case* 1. $F(A; X) = |X_1| + |X_2| + \cdots + |X_n|$. In this case, we have

$$(4.13) \qquad \sum_{i=1}^{n} |X_i| \quad \sim \quad \sum_{j=1}^{n} |\sigma_j|^{1/j}.$$

A proof may be found in [18].

*Case* 2. Suppose $F$ is of the form $F(A; X) = \sum_i |f_i(A; X)|^\varepsilon$. We can write $f_i(A; X) = F_i(X)$ for some $F_i \in L[X]$ where $L = C[A]$. Then

$$F(A; X) = \frac{1}{n!} \sum_{\pi \in S_n} F(A; X^\pi) = \frac{1}{n!} \sum_\pi \sum_i |F_i(X^\pi)|^\varepsilon \sim (\sum_i \sum_\pi |F_i(X^\pi)|)^\varepsilon.$$

Fix an ordering of $S_n$: $S_n = \{\pi_1, \ldots, \pi_{n!}\}$. Then, applying (4.13), we have, for each fixed $i$:

$$\sum_\pi |F_i(X^\pi)| = \sum_{j=1}^{n!} |F_i(X^{\pi_j})| \sim \sum_{j=1}^{n!} |\sigma_j^{(n!)}|^{1/j}$$

where

$$\sigma_j^{(n!)} = \sigma_j^{(n!)}(F_i(X^{\pi_1}), \ldots, F_i(X^{\pi_{n!}})) = \Phi(X_1, \ldots, X_n)$$

with $\Phi \in L[X_1, \ldots, X_n]$. It is clear that $\Phi(X)$ is symmetric. Hence, by Lemma 4.2, it is an element of $L[\sigma_1, \ldots, \sigma_n]$ and this completes the proof of Case 2.

*Remark.* Let $F \in C|(A; X)|$ be a general ARP which is symmetric in the $X$ variables. Then we claim that the conclusion of Lemma 4.3 holds: We do not know that the numerator and denominator of $F$ are each symmetric (in fact they need not be) so we cannot immediately apply Lemma 4.3. But since $F(A; X) = \frac{1}{n!} \sum_{\pi \in S_n} F(A; X^\pi)$ we see that we can re-write $F$ in a form where the denominator is symmetric, and the numerator is symmetric, and thus we can apply the lemma to the numerator and denominator separately.

*Proof of Lemma* 4.1. Using the observation (4.2a), we see that we may assume that $\varepsilon = \delta = 1$. Now let

$$R(Z; B) = \frac{\sum_{i=1}^{I} |\sum_{\mu=0}^{M_i} B'_{i,\mu} Z^\mu|}{\sum_{j=1}^{J} |\sum_{\nu=0}^{N_j} B''_{j,\nu} Z^\nu|}$$

and let $S(X_1, \ldots, X_d; Z) = \sum_{i=1}^{d} R(X_i; B)$. Then, if we express $S$ as a single fraction whose denominator is the product of the denominators of the $R(X_i; B)$, we see that the numerator and denominator of $S$ are each symmetric in the $X$ variables and so, applying Lemma 4.2 twice, once in the numerator and once in the denominator, we write $S(X; B) = \tilde{S}(\sigma_1, \ldots \sigma_d; B)$. Now we define $T(A; B)$ as follows: $T(A; Z) = \tilde{S}(A_1/A_0, \ldots, A_d/A_0; B)$. One easily sees that $T$ satisfies (4.12) and thus Lemma 4.1 is proved.

C. *The case of a single polynomial in the denominator with simple roots.* We begin with a preliminary theorem which treats the simplest case: that in which the denominator is a scalar polynomial with distinct roots. We shall continue using the notation established in Section 3.A:

Fix $\varepsilon, \delta$ nonnegative rational numbers and assume $\delta < 2$, $I = 1$, $J = 1$. Let $M = M_1$ and $N = N_1$, and let $\mathcal{B} = \{b = (b', b'')\}$ be the coefficient space defined by (4.4). To each $b \in \mathcal{B}$ we associate $R_b(Z) = |P(b; Z)|^{\varepsilon} |Q(b; Z)|^{-\delta} \in \mathbb{C}|(Z)|$ as in (4.5). Assume $Q(b; z)$ has distinct roots. Assume further that $(\varepsilon, \delta)$ is nondegenerate; that is, $\nu\varepsilon - (N - k)\delta + 2 \neq 0$ for all $\nu, k$ in the ranges $0 \leq \nu \leq M$ and $0 \leq k \leq N - 1$.

THEOREM 5. *There exists $T(B) = T(B', B'') \in \mathbb{Z}|(B', B'')|$ and $s_1 \in (0, 1)$ such that whenever $|||Q(b; Z) - Z^N||| < s_1$, the following holds:*

1. *The vector $b$ is in the strict domain of $T$.*

2. *The following estimate holds:*

$$(4.14) \qquad \int_B R_b(z) \ dV \ = \ \int_B \frac{|P(b, z)|^{\varepsilon}}{|Q(b, z)|^{\delta}} \ dV \quad \sim \quad T(b', b'') \ = \ T(b) \ ,$$

*where the implied constants are independent of $b$.*

*Proof.* By Rouché's theorem, there exists an $s_1 \in (0, 1)$ such that $|||Q(Z) - Z^N||| < s_1$ implies that all the roots of $Q(Z)$ are of size less than $1/2$. Fix such an $s_1$. We also made the assumption that all the roots of $Q$ are simple. Thus the version (2.20) of Theorem 1 applies. We can now recast (2.20) as an ARP in the coefficients of $P(z)$ and $Q(z)$, using arguments similar to those of the proof of (i), Theorem 3. We shall do this in a systematic fashion, using Lemmas 4.1-4.3: Define, for $0 \leq k \leq N - 2$, the following absolute polynomial power:

$$\tilde{F}_k(Y; X_1, \ldots, X_N) = \sum_{\lambda \in \Lambda_k} \prod_{\nu=0}^{k} |Y - X_{i_\nu}|$$

where $\Lambda_k = \{\lambda = (i_0, \ldots, i_k) : 1 \leq i_\nu \leq N, \quad i_\nu \neq i_\mu \ \text{if} \ \mu \neq \nu\}$. Set $\tilde{F}_k(Y; X_1, \ldots, X_N) = 1$ if $k < 0$.

We know, by Lemma 4.3, that $\tilde{F}_k(Y; X_1, \ldots, X_N) = F_k(Y; \Sigma)$, where $\Sigma = \Sigma(X)$ is the vector of elementary symmetric polynomials in the $X$ variables. Thus, in view of (3.13) and Lemma 4.1, we obtain

$$L_0(\alpha) L_1(\alpha) \cdots L_k(\alpha) \ \sim \tilde{F}_k(\alpha; \alpha_1, \ldots, \alpha_N) \sim F_k(\alpha; b''/b_N'')$$

where we have divided $Q(z)$ by $n''_N$ in order to make it monic (note that $b''_N \neq 0$ for $s_1 < 1$). Now use (4.2a) to choose $H(Y; B) = H(Y; B'; B'') \in \mathbb{Z}|(Y; B'; B'')|$ with the following property:

$$H(Y; B'; B'') \sim \sum_{\nu=0}^{M} |P^{(\nu)}(Y)|^{\varepsilon} \sum_{k=-1}^{N-2} \frac{F_k(Y; B''/B''_N)^{\nu\varepsilon-(N-k)\delta+2}}{F_{k-1}(Y; B''/B''_N)^{\nu\varepsilon-(N-k-1)\delta+2}}$$

where $P(Y) = P(Y; B') = \sum_{\mu=0}^{M} B'_\mu Y^\mu$. The assumption that $Q$ has no multiple roots implies that $(\alpha, b', b'')$ is in the strict domain of $H$.

Next let $\tilde{T}_1(Y_1, \ldots, Y_N; B'; B'') = \sum_{i=1}^{N} H(Y_i; B'; B'')$. Then, by Lemma 4.3, we know that $\tilde{T}_1(Y; B'; B'') = T_1(\Sigma(Y); B'; B'')$. Therefore we have

$$\sum_{\{\alpha: Q(\alpha)=0\}} H(\alpha; b', b'') \sim T_1(b''/b''_N; b'; b'').$$

Define $T \in \mathbb{Z}|(B'; B'')|$ by $T(B'; B'') = T_1(B''/B''_N; B'; B'')$. Now set $\Lambda = 1$ in (2.6). Since $Q(z)$ has no multiple roots, we can write for all $(\nu, k)$ with $k \leq N-2$

$$\frac{1}{\Phi_{\nu,k}(\alpha)} = \frac{(\prod_{i=0}^{k} L_i(\alpha))^{\nu\varepsilon-(N-k)\delta+2}}{(\prod_{i=0}^{k-1} L_i(\alpha))^{\nu\varepsilon-(N-k-1)\delta+2}}.$$

Theorem 5 follows.

D. *The regularization process for multiple roots.* We now come to one of the most important features of Theorem 4, namely the emergence of a stratification in the general case. This stratification is a reflection of the interplay between $P(z)$ and $Q(z)$ near the roots of $Q(z)$. At the technical level, it is a consequence of the fact that the estimate (2.20) breaks down when $Q(z)$ has multiple roots. Equivalently, if we consider the case of multiple roots as limiting cases of simple roots, then the ARP $T(B', B'') \in \mathbb{Z}|(B', B'')|$, describing the integral of $|P(z)|^{\varepsilon}|Q(z)|^{-\delta}$, will approach $0/0$. To obtain a formula similar to (2.20) in the case of multiple roots, which is also valid for a general ARP, we proceed in three steps: The first is a certain "regularization process," which is justified by the Lebesgue Monotone Convergence Theorem:

$$(4.15) \qquad \int_B \frac{\sum_{i=1}^{I} |P_i(z)|^{\varepsilon}}{\sum_{j=1}^{J} |Q_j(z)|^{\delta}} dV = \lim_{\mu \downarrow 0} \int_B \frac{\sum_{i=1}^{I} |P_i(z)|^{\varepsilon}}{(\sum_{j=1}^{J} |Q_j(z)| + \mu)^{\delta}} dV.$$

Here $\mu$ is a positive parameter which approaches zero from the right.

The next step is the introduction of "theta parameters," which will allow us to reduce to the case where there is just a single term in the denominator, that is, the case $J = 1$:

$$(4.16) \qquad \int_B \frac{\sum_{i=1}^I |P_i(z)|^\varepsilon}{(\sum_{j=1}^J |Q_j(z)| + \mu)^\delta} dV$$

$$\sim \int_{(\mathbb{R}/\mathbb{Z})^J} d\theta_1 \cdots d\theta_J \int_B \frac{\sum_{i=1}^I |P_i(z)|^\varepsilon}{|\sum_{j=1}^J Q_j(z) e^{2\pi i \theta_j} + \mu|^\delta} dV.$$

This is justified, for $\delta < 1$, by the following elementary lemma:

LEMMA 4.4.     *Let $a_1, \cdots, a_J$ be a finite set of complex numbers, and $\delta$ a real number with $0 < \delta < 1$. Then*

$$(4.17)$$
$$\int_{(\mathbb{R}/\mathbb{Z})^J} \frac{d\theta_1 \cdots d\theta_J}{|\sum_{j=1}^J a_j e^{2\pi i \theta_j}|^\delta} = \int_{(\mathbb{R}/\mathbb{Z})^{J-1}} \frac{d\theta_1 \cdots d\theta_{J-1}}{|\sum_{j=1}^{J-1} a_j e^{2\pi i \theta_j} + a_J|^\delta} \sim \frac{1}{(\sum_{j=1}^J |a_j|)^\delta},$$

*where the implied constants depend only on $\delta$ and $J$.*

The third step consists of sampling lemmas, which allow us to calculate the size of the integral (4.14) by sampling a finite number of theta parameters; a discussion of this step is postponed until Section F.

*Proof of Lemma* 4.4. The case $J = 1$ is trivial, so we start with $J = 2$. We must show

$$\int_{\mathbb{R}/\mathbb{Z}} \frac{1}{|a + be^{2\pi i \theta}|^\delta} d\theta \sim \frac{1}{|a|^\delta + |b|^\delta}.$$

Without loss of generality, we may assume that $a$ is a positive real number, that $b = 1$ and $a \leq b$. However, in this case, we have

$$\int_0^1 \frac{d\theta}{2^\delta} \leq \int_0^1 \frac{d\theta}{|a + e^{2\pi i \theta}|^\delta} \leq \int_0^1 \frac{d\theta}{|\sin 2\pi \theta|^\delta} < \infty.$$

Turning to the general case, we may reduce as before to the case where $0 < a_j \leq 1$ and $a_J = 1$. Integrating first with respect to $\theta_J$, using the $J = 2$ result, we obtain

$$\int_{(\mathbb{R}/\mathbb{Z})^J} \frac{d\theta_1 \cdots d\theta_J}{|\sum_{j=1}^J a_j e^{2\pi i \theta_j}|^\delta} \sim \int_{(\mathbb{R}/\mathbb{Z})^{J-1}} \frac{d\theta_1 \cdots d\theta_{J-1}}{|\sum_{j=1}^{J-1} a_j e^{2\pi i \theta_j}|^\delta + 1} \sim 1.$$

The lemma is proved.

In order to apply Theorem 5, we need to be assured that the denominator in (4.16) has only simple roots. This will be guaranteed by the following:

LEMMA 4.5.     *Let $Q(Z) \in \mathbb{C}[Z]$ and for $\mu > 0$, let $Q_\mu(Z) = Q(Z) + \mu$. For $\mu$ sufficiently small and different from zero, $Q_\mu(z)$ has only simple roots.*

*Proof.* The discriminant of the polynomial $Q(Z)$ is a polynomial in $\mu$ which does not vanish identically. Hence, it can vanish only for a finite number of values of $\mu$.

E. *The case of a single polynomial in the denominator.* In this section we prove Theorem 4 in the case where $I = J = 1$. By the first remark following the statement of Theorem 4, we may assume that $\delta < 1$. For $\mu > 0$ and $b \in \mathcal{B}$ let

$$(4.18) \quad \mathcal{I}_b = \int_B \frac{|P(b;z)|^\varepsilon}{|Q(b;z)|^\delta} \, dV \quad \text{and} \quad \sigma_b(\mu) = \int_B \frac{|P(b;z)|^\varepsilon}{(|Q(b;z)| + \mu)^\delta} \, dV.$$

As observed in the preceding section:

$$(4.19) \quad \mathcal{I}_b = \lim_{\mu \downarrow 0} \sigma_b(\mu) \ .$$

Thus our task is to evaluate the above limit of $\sigma_b(\mu)$. Let

$$(4.20) \quad E_b(\theta, \mu) = \int_B \frac{|P(b;z)|^\varepsilon}{|Q(b;z) + \mu e^{2\pi i \theta}|^\delta} dV.$$

Then, in view of Lemma 4.4 (and (4.16)), we may write

$$(4.21) \quad \sigma_b(\mu) \sim \int_{\mathbb{R}/\mathbb{Z}} E_b(\theta, \mu) d\theta,$$

where here, as is the case with all equivalences in this section, the implied constants are independent of $\mu$ and $b$. It is easy to see that

$$(4.22) \quad \sigma_b(\mu) \quad \sim \quad \inf_{\theta \in \mathbb{R}/\mathbb{Z}} E_b(\theta, \mu).$$

In fact, it is evident that $\sigma_b(\mu) \leq \inf_\theta E_b(\theta, \mu)$. On the other hand, if $\inf_\theta E_b(\theta, \mu) \geq K\sigma_b(\mu)$, it would follow that

$$\int_{\mathbb{R}/\mathbb{Z}} E_b(\theta, \mu) d\theta \geq K\sigma_b(\mu),$$

which would contradict the estimate (4.21), if $K$ is large compared to the universal constants implicit in that estimate. This establishes (4.22).

Let $w = \mu e^{2\pi i \theta}$ and let $b(w) = (b'_M, \dots b'_0, b''_N, \dots, b''_0 + w)$. Then we have the obvious identity: $Q(b;z) + w = Q(b(w);z)$. With this notation, we can write

$$E_b(\theta, \mu) = \int_B \frac{|P(b;z)|^\varepsilon}{|Q(b(w);z)|^\delta} dV.$$

Now assume $b \in \mathcal{B}_0$ where $\mathcal{B}_0$ is the (Euclidean) open subset of $\mathcal{B}$ defined by

$$(4.22a) \quad \mathcal{B}_0 \quad = \quad \{b \in \mathcal{B} : |||Q(b;Z) - Z^N||| < s_1/2\} \ .$$

Here $s_1$ is the parameter which appears in Theorem 5. Assume $0 < \mu < s_1/2$. Then $|||Q(b(w);Z) - Z^N||| < s_1$ . Moreover, for $\mu$ sufficiently small, we know, by Lemma 4.5, that $Q(b(w); z)$, has only simple roots. Thus, applying Theorem 5 we obtain, for $b \in \mathcal{B}_0$ and $\mu$ sufficiently small:

$$(4.23) \qquad E_b(\theta, \mu) \sim T(b(w)) \; = \; \tilde{T}(b, w)$$

where $\tilde{T}(B; W) = \tilde{T}(B', B''; W) \in \mathbb{Z}|(B; W)|$ is defined as follows:

$$(4.24) \quad \tilde{T}(B', B'', W) = T(B'_M, B'_{M-1}, \ldots, B'_0, B''_N, B''_{N-1}, \ldots, B''_1, B''_0 + W)$$

and $T(B)$ is as in Theorem 5. We note that in (4.23), the implied constants are independent of $b, \mu$ and $\theta$.

Estimates (4.22) and (4.23) imply that for $b \in \mathcal{B}_0$ and $\mu$ sufficiently small,

$$(4.25) \qquad \sigma_b(\mu) \sim \inf_{\theta \in \mathbb{R}/\mathbb{Z}} \tilde{T}(b; \mu e^{2\pi i \theta})$$

and thus, letting $\mu$ tend to 0,

$$\mathcal{I}_b \; = \; \lim_{\mu \downarrow 0} \sigma_b(\mu) \; = \; \lim_{\mu \downarrow 0} \inf_{\theta \in \mathbb{R}/\mathbb{Z}} \tilde{T}(b; \mu e^{2\pi i \theta}).$$

On the other hand, by observation 2 of Section 4.A, $\lim_{w \to 0} \tilde{T}(b; w)$ exists. Thus, for all $b \in \mathcal{B}_0$ we have

$$(4.26) \qquad \mathcal{I}_b = \lim_{w \to 0} \tilde{T}(b; w) \; ,$$

where the implied constants are independent of $b$.

It remains to clarify the dependence of the limit $\lim_{w \to 0} T(b; w)$ on $b$, the vector of coefficients of the polynomials $P(b; Z)$ and $Q(b; Z)$. In order to do so, we write the size of $T(B; W)$ in the form

$$(4.27) \qquad T(B; W) \; \sim \; \left( \frac{\sum_{i=1}^{I^*} |\sum_{m=0}^{M_i} F_{im}(B) W^m|}{\sum_{j=1}^{J^*} |\sum_{n=0}^{N_j} G_{jn}(B) W^n|} \right)^{\varepsilon_1},$$

where the coefficients $F_{im}(B)$ and $G_{jn}(B)$ are polynomials in $\mathbb{C}[B]$ . Let $m_0$ be the smallest integer $m$ such that $F_{im} \neq 0$ for some $i$. For $\lambda$ a nonnegative integer, define:

$$\tilde{U}_\lambda \; = \{ b \in \mathcal{B}_0 : \sum_{i=1}^{I^*} \sum_{m=0}^{m_0 + \lambda - 1} |F_{im}(b)| \; = \; 0 \}.$$

Let $U_\lambda$ be the Zariski closure of $\tilde{U}_\lambda$ in $\mathcal{B}$ (that is, the smallest variety containing $\tilde{U}_\lambda$). Then the $U_\lambda$ form a finite decreasing chain: $\mathcal{B} = U_0 \supseteq U_1 \supseteq \cdots \supseteq U_\mathcal{N}$.

The fact that $\mathcal{I}_b > 0$, together with (4.26) implies that $U_\mathcal{N} = \emptyset$.

Let $b \in \mathcal{B}_0$. Then there exists a unique $\lambda$ such that $b \in U_\lambda \backslash U_{\lambda+1}$. Again, (4.26) and the positivity of $\mathcal{I}_b$ imply

$$(4.28) \qquad \sum_{i=1}^{J^*} \sum_{m=0}^{m_0 + \lambda - 1} |G_{im}(b)| \; = \; 0;$$

thus we define

$$(4.29) \qquad T_\lambda(B) \quad = \quad \left( \frac{\sum_{i=1}^{I^*} \sum_{m=0}^{m_0+\lambda} |F_{im}(B)|}{\sum_{j=1}^{J^*} \sum_{m=0}^{m_0+\lambda} |G_{jm}(B)|} \right)^{\varepsilon_1}.$$

Combining (4.26), (4.27), (4.28) and (4.29), we conclude that for $b \in \mathcal{B}_0$ and $b \in U_\lambda \backslash U_{\lambda+1}$, we have

$$(4.30) \qquad\qquad\qquad \mathcal{I}_b \quad \sim \quad T_\lambda(b).$$

This proves the first two parts of Theorem 4. The third part follows from the construction of the $U_\lambda$.

*Example.* It may be easier to understand the stratification of Theorem 4 with the following easier example in mind. Consider the integral

$$(4.31) \qquad\qquad \mathcal{I}(a,b,c) = \int_{B_\Lambda} \frac{|z-c|^\varepsilon}{|az^2 - bz|^\delta} dV$$

with $1 < \delta$, and $2\delta - \varepsilon < 2$ and $(a,b,c) \neq (0,0,0)$. Then $\mathcal{I}(a,b,c) = \infty$ if and only if $b = 0$ and $c \neq 0$. Moreover, the size of the integral can be given as follows:

$$(4.32) \qquad \mathcal{I}(a,b,c) \sim (\Lambda|a| + |b|)^{-\delta} \Lambda^{-\delta+2+\varepsilon} + \left( \frac{\Lambda}{\Lambda|a| + |b|} \right)^{2-\delta} |b|^{-2\delta+2} |c|^\varepsilon$$

for $b \neq 0$. For $b = 0$, we have $\mathcal{I}(a,b,c) = \infty$ if $c \neq 0$, and for $b = c = 0$, we have $\mathcal{I}(a,b,c) \sim |a|^{-\delta} \Lambda^{-2\delta+2+\varepsilon}$.

*Remarks.* 1. In this paper, we have chosen to use $\theta$ parameters, as these parameters are also instrumental in the treatment of the case of several polynomials in the denominator (see the next section). But we can use also the following formula for the regularized integral (4.18) which may be of independent interest

$$(4.33) \quad \int_{B_\Lambda} \frac{|P(z)|^\varepsilon}{(|Q(z)| + \mu^N)^\delta} dV$$

$$\sim \sum_{\{\alpha; Q(\alpha)=0\}} \sum_{\nu=0}^{M} |P^{(\nu)}(\alpha)|^\varepsilon \prod_{k=0}^{N-1} \left[ \mu^N + L_k(\alpha)^{N-k} \prod_{0 \leq i < k} L_i(\alpha) \right]^{-\kappa_k}.$$

Here $P(z)$ and $Q(z)$ are polynomials of degrees $M$ and $N$ respectively, $Q(z)$ is monic, the exponents $\kappa_k$ are defined by $\kappa_k = \frac{\nu\varepsilon+2}{(N-k-1)(N-k)}$ for $0 \leq k \leq N-2$, $\kappa_{N-1} = \delta - \nu$, and we have assumed $M\nu + 2 < \delta N$ for simplicity.

2. The function $T_\lambda(b)$ is infinite on the variety $\{b \in \mathcal{U}_\lambda \backslash \mathcal{U}_{\lambda+1}; L_\lambda(b) = 0\}$. By intercalating such varieties, we could also have formulated Theorem 4 in terms of a seemingly finer stratification $\mathcal{U}_\mu^*$, where the values of the ARP on $\mathcal{U}_\mu^* \backslash \mathcal{U}_{\mu+1}^*$ would be either finite or infinite throughout.

F. *The sampling lemmas.* In order to treat general denominators in ARP's, and not just those consisting of a single absolute value, we first introduce theta parameters, which reduce a denominator of the form $\sum_{j=1}^{J} |Q_j(z)|$ to one of the form $|\sum_{j=1}^{J} Q_j(z)e^{2\pi i\theta_j}|$ for $\theta_j \in \mathbb{R}/\mathbb{Z}$. We then apply a certain *sampling lemma*, which replaces the integral over $\mathbb{R}/\mathbb{Z}$ to a sum over a boundedly finite number of $\theta_j$. We begin by stating such a lemma in its simplest form.

LEMMA 4.6 (first sampling lemma).    *Fix nonnegative rational numbers, $\varepsilon$, $\delta$, such that $\delta < 1$ and $(\varepsilon, \delta)$ is nondegenerate. Fix nonnegative integers $M$ and $N$. Then there exists a real number $s \in (0, 1)$, and a positive integer $d$, depending only on $\varepsilon, \delta, M, N$, with the following property*: *For all polynomials $P(z)$, $Q_1(z)$, $Q_2(z)$ with*

- $\deg(P) = M$, $\deg(Q_2) \leq \deg(Q_1) = N$;
- $|||Q_1(Z) - Z^N||| < s/2$ *and* $|||Q_2(Z)||| < s/2$,

*we have*

$$(4.34)$$
$$\int_B \frac{|P(z)|^\varepsilon}{(|Q_1(z)| + |Q_2(z)|)^\delta} dV \;\sim\; \inf_{\{\theta \in \mathbb{R}/\mathbb{Z}\}} \int_B \frac{|P(z)|^\varepsilon}{|Q_1(z) + Q_2(z)e^{2\pi i\theta}|^\delta} dV$$

$$\sim\; \inf_{\{\theta \in \mathbb{R}/\mathbb{Z}: d\theta = 0\}} \int_B \frac{|P(z)|^\varepsilon}{|Q_1(z) + Q_2(z)e^{2\pi i\theta}|^\delta} \; dV.$$

*Here the implied constants depend only on the fixed constants $\varepsilon, \delta, M, N$, and $\{\theta \in \mathbb{R}/\mathbb{Z} : d\theta = 0\}$ is the set of points of the form $t = \frac{k}{d}$, $0 \leq k \leq d-1$.*

*Proof.* In the proof, all equivalences will have implied constants which depend only on the fixed constants $\varepsilon, \delta, M, N$. The lemma is obviously true in the case when $P(z)$ is the zero polynomial, and all expressions in the estimate (4.34) vanish. Thus we assume that $P(z)$ is not identically zero, and that all the terms figuring in (4.34) are strictly positive. We clearly have

$$(4.35)$$
$$\int_B \frac{|P(z)|^\varepsilon}{(|Q_1(z)| + |Q_2(z)|)^\delta} dV \;\leq\; \inf_{\{\theta \in \mathbb{R}/\mathbb{Z}\}} \int_B \frac{|P(z)|^\varepsilon}{|Q_1(z) + Q_2(z)e^{2\pi i\theta}|^\delta} dV$$

$$\leq\; \inf_{\{\theta \in \mathbb{R}/\mathbb{Z}: d\theta = 0\}} \int_B \frac{|P(z)|^\varepsilon}{|Q_1(z) + Q_2(z)e^{2\pi i\theta}|^\delta} \; dV,$$

so our main task is to show that

$$(4.36)$$
$$\inf_{\{\theta \in \mathbb{R}/\mathbb{Z}: d\theta = 0\}} \int_B \frac{|P(z)|^\varepsilon}{|Q_1(z) + Q_2(z)e^{2\pi i\theta}|^\delta} \; dV \leq C \int_B \frac{|P(z)|^\varepsilon}{(|Q_1(z)| + |Q_2(z)|)^\delta} dV$$

for a constant $C$ depending only on $\varepsilon$, $\delta$, $M$, $N$. Let $s \in (0,1)$ be chosen as in Theorem 4, in the case we have already proven, namely when $|I| = |J| = 1$. Let $P(z), Q_1(z)$ and $Q_2(z)$ be polynomials satisfying the conditions in the lemma, and consider, for $w \in \mathbb{C}$, the following expression

$$J(w) = \int_B \frac{|P(z)|^\varepsilon}{|Q_1(z) + Q_2(z)w|^\delta}\, dV.$$

Thus $J(w) = \int_B R_{b(w)}dV$ where $b : \mathbb{C} \to \mathcal{B}$ is a complex line, $L$, in our parameter space $\mathcal{B}$.

Let $\mathcal{B} = U_0 \supset U_1 \supseteq \cdots \supset U_\mathcal{N} = \emptyset$ be the filtration defined in section E. Then there exists $\lambda$ such that $L \subseteq U_\lambda$ but $L \not\subseteq U_{\lambda+1}$. Since $L$ is one-dimensional, and since each $U_{\lambda+1}$ is defined by a finite set of algebraic equations, $L \cap U_{\lambda+1}$ consists of a finite set of points. Furthermore, if $w \in \mathbb{C}$ has length $|w| \leq 1$, then $Q(b(w), Z) = Q_1(b_1, Z) + wQ_2(b_2, Z)$ has degree $N$ and $|||Q(Z) - Z^{N_1}||| < s$. Thus the case $|I| = |J| = 1$ of Theorem 4 applies, and we see that there is a finite set $S \subseteq \mathbb{C}$ such that

(4.37)
$$\int_B \frac{|P(z)|^\varepsilon}{|Q_1(z) + wQ_2(z)|^\delta}\, dV \sim T_\lambda(b(w)) \sim H(w) \quad \text{for all} \quad w \notin S, \quad |w| \leq 1,$$

where, either $H(w) \equiv \infty$, or $H(W) \in \mathbb{C}|(W)|$ is a one variable ARP. We may assume that $H$ has the form given by (4.2a) . Note that while $H$ does, in general, depend on $b$, the number of terms which appear in $H$, as well as the degrees of those terms, is bounded by a constant $D$ which depends only on $\varepsilon, \delta, M$ and $N$.

As we have seen in Section A, $H(w)$ admits a continuous extension, possibly $\infty$ valued, to all of $|w| \leq 1$. However, the left-hand side of (4.37) is not continuous in $w$, and we must proceed with some care. Since an arbitrary value $w$ in $|w| \leq 1$ can be approached by values outside of $S$, Fatou's lemma, and the fact that $H(w)$ is continuous, implies that there exists $C > 0$, depending only on $\varepsilon, \delta, M$ and $N$, such that

(4.38)
$$\int_B \frac{|P(z)|^\varepsilon}{|Q_1(z) + wQ_2(z)|^\delta}\, dV \leq C\, H(w), \quad \text{for all } w, \ |w| \leq 1.$$

We restrict now our attention to $|w| = 1$, and write $w = e^{2\pi i\theta}$. As in the case considered in Section 3.E, it is convenient to introduce the notation

(4.39)
$$\sigma = \int_B \frac{|P(z)|^\varepsilon}{(|Q_1(z)| + |Q_2(z)|)^\delta}\, dV$$

$$E(\theta) = \int_B \frac{|P(z)|^\varepsilon}{|Q_1(z) + e^{2\pi i\theta}Q_2(z)|^\delta}\, dV.$$

We observe that

$$(4.40) \qquad \sigma \sim \inf_{(\mathbb{R}/\mathbb{Z})\setminus S} E(\theta),$$

where the implied constant depends only on $\delta$. The argument is the same as in Section E. Clearly, $\sigma \leq E(\theta)$ for all $\theta$, and in particular $\sigma \leq \inf_{(\mathbb{R}/\mathbb{Z})\setminus S} E(\theta)$. On the other hand, by Lemma 4.4 and the fact that $S$ is of measure 0, we see that

$$\sigma \sim \int_{\mathbb{R}/\mathbb{Z}} E(\theta) d\theta = \int_{(\mathbb{R}/\mathbb{Z})\setminus S} E(\theta) d\theta \geq \inf_{(\mathbb{R}/\mathbb{Z})\setminus S} E(\theta).$$

This establishes (4.40). Together with the estimate (4.37) and the continuity of $H(w)$, this implies

$$(4.41) \qquad \sigma \sim \inf_{(\mathbb{R}/\mathbb{Z})\setminus S} H(e^{2\pi i\theta}) = \inf_{\mathbb{R}/\mathbb{Z}} H(e^{2\pi i\theta}).$$

The key fact which we require now is that the infimum of $H(e^{2\pi i\theta})$ over the full circle $\theta \in \mathbb{R}/\mathbb{Z}$ can actually be determined by its infimum over a boundedly finite set division points, $\{\theta \in \mathbb{R}/\mathbb{Z} : d\theta = 0\}$. The precise statement and proof will be given in Lemma 4.8. First, we establish a preliminary lemma, which says that $H$ has the same size as its infimum on a set of large measure:

LEMMA 4.7. *There exists a constant $K$ depending only on $\varepsilon, \delta, M, N$ such that*

$$(4.42) \qquad \mathrm{Vol}\{\theta \in \mathbb{R}/\mathbb{Z};\ H(e^{2\pi i\theta}) \leq K \inf_{\mathbb{R}/\mathbb{Z}} H(e^{2\pi i\theta})\} \geq 1/2,$$

*where* Vol *indicates the Euclidean measure.*

*Proof of Lemma* 4.7. The argument is similar to the one for $E(\theta)$ in (4.22). On one hand, we have

$$\int_{\mathbb{R}/\mathbb{Z}} H(e^{2\pi i\theta}) d\theta = \int_{(\mathbb{R}/\mathbb{Z})\setminus S} H(e^{2\pi i\theta}) d\theta \sim \int_{(\mathbb{R}/\mathbb{Z})\setminus S} E(\theta)\, d\theta = \int_{\mathbb{R}/\mathbb{Z}} E(\theta) d\theta \sim \sigma.$$

On the other hand,

$$\int_{\mathbb{R}/\mathbb{Z}} H(e^{2\pi i\theta}) d\theta \geq K\sigma\, \mathrm{Vol}\{\theta \in \mathbb{R}/\mathbb{Z} : H(e^{2\pi i\theta}) > K\sigma\}.$$

For $\sigma$ strictly positive and finite, and for $K$ large compared to the implicit constant in the earlier inequality (which depends only on $\varepsilon, \delta, M$ and $N$), we obtain

$$\mathrm{Vol}\{\theta; H(e^{2\pi i\theta}) \leq K\sigma\} \geq 1/2.$$

Since $\sigma \sim \inf_{\mathbb{R}/\mathbb{Z}} H(e^{2\pi i\theta})$ in view of (4.41), the lemma is proved in this case. Evidently, $\sigma$ cannot vanish unless $P(Z) \equiv 0$, and if $\sigma = \infty$, then the estimate (4.42) is trivial. Thus Lemma 4.7 is proved. We now establish a sampling lemma for ARP's, which, when combined with (4.38) and (4.41), will complete the proof of Lemma 4.6:

LEMMA 4.8.    *Fix a positive number $\varepsilon^{\#}$, an integer $d^{\#}$, and a positive constant $K > 1$. Then there exists a positive integer $d$ depending only on $d^{\#}$, such that for all $H(W) \in \mathbb{C}|(W)|$ of the form*

$$(4.43) \qquad H(W) = \left( \frac{\sum_{i=1}^{I} |P_i(W)|^2}{\sum_{j=1}^{J} |Q_j(W)|^2} \right)^{\varepsilon^{\#}/2},$$

*with $|I|, |J|$ and the degrees of $P_i(W)$ and $Q_j(W)$ all bounded by $d^{\#}$, and satisfying the following two conditions:*

- $\inf_{\mathbb{R}/\mathbb{Z}} H(e^{2\pi i\theta}) > 0$;
- $\mathrm{Vol}\{\theta \in \mathbb{R}/\mathbb{Z} : H(e^{2\pi i\theta}) \leq K \inf_{\mathbb{R}/\mathbb{Z}} H(e^{2\pi i\theta})\} \geq 1/2$,

*we have*

$$(4.44) \qquad \inf_{\mathbb{R}/\mathbb{Z}} H(e^{2\pi i\theta}) \sim \inf_{\{\theta \in \mathbb{R}/\mathbb{Z} : d\theta = 0\}} H(e^{2\pi i\theta}).$$

*with the implied constants depending only on $d^{\#}$ and $K$.*

*Proof of Lemma 4.8.* Let $\tau = \inf_{\mathbb{R}/\mathbb{Z}} H(e^{2\pi i\theta})$. We may assume that $\tau$ is finite; otherwise the lemma is trivial. We claim that the set

$$(4.45) \qquad \{\theta \in \mathbb{R}/\mathbb{Z} : H(e^{2\pi i\theta}) \leq K\tau\}$$

consists of a finite number of intervals, and that the number of intervals is bounded by a constant depending only on $d^{\#}$: It suffices to prove that $\{\theta : H(e^{2\pi i\theta}) = K\tau\}$ has cardinality bounded uniformly in terms of $d^{\#}$. To see this, we rewrite the equation for these points as

$$(4.46) \qquad \sum_{i=1}^{I} |P_i(e^{2\pi i\theta})|^2 = (K\tau)^{(2/\varepsilon^{\#})} \sum_{j=1}^{J} |Q_j(e^{2\pi i\theta})|^2.$$

If we multiply both sides of the equation by $e^{2\pi i d^{\#}\theta}$, each of the terms

$$e^{2\pi i d^{\#}\theta} |P_i(e^{2\pi i\theta})|^2 = e^{2\pi i d^{\#}\theta} P_i(e^{2\pi i\theta}) \overline{P}_i(e^{-2\pi i\theta})$$

and

$$e^{2\pi i d^{\#}\theta} |Q_j(e^{2\pi i\theta})|^2 = e^{2\pi i d^{\#}\theta} Q_j(e^{2\pi i\theta}) \overline{Q}_j(e^{-2\pi i\theta})$$

becomes a polynomial in the variable $w = e^{2\pi i\theta}$. Thus the equation (4.46) is a polynomial equation for $w$ of degree at most $2d^{\#}$. It cannot have infinitely many solutions, since then two sides would then be identically equal, and we would have $H(e^{2\pi i\theta}) \equiv K\tau$. For $K > 1$, this contradicts the fact that $\tau$ is the infimum of $H(e^{2\pi i\theta})$. Thus the number of solutions of (4.46) is finite, and bounded by $2d^{\#}$. It follows that the number of disjoint intervals making up the set (4.45) is bounded by $2d^{\#}$.

Consider now the set of points $\theta$ satisfying $d\theta = 0$. If this set does not intersect the set (4.45), it would imply that each of the intervals making up the set (4.45) has length at most $1/d$. But we just saw that the set (4.45) consists of at most $d^{\#}$ intervals. Thus its total measure would be bounded by $d^{\#}/d$. Hence, if we choose $d$ large enough so that $(d^{\#}/d) < (1/2)$, we obtain a contradiction. Thus, we conclude that for $d > 2d^{\#}$, there exists $\theta \in \mathbb{R}/\mathbb{Z}$ such that $d\theta = 0$ and such that $T(e^{2\pi i\theta}) \leq K\tau$. Lemma 4.8 is proved.

We return now to the proof of the sampling lemma. Recall that we may assume without loss of generality that $P(z)$ is not identically zero. In this case the integral $\sigma$ of (4.39) is strictly positive; hence the infimum of $T(e^{2\pi i\theta})$ is strictly positive in view of (4.41). In view of Lemma 4.7, the function $T(w)$ satisfies all requirements in Lemma 4.8. Thus we may apply (4.41) and Lemma 4.8 and deduce that

$$(4.47) \qquad \sigma \sim \inf_{\mathbb{R}/\mathbb{Z}} H(e^{2\pi i\theta}) \sim \inf_{\{\theta \in \mathbb{R}/\mathbb{Z};\, d\theta=0\}} H(e^{2\pi i\theta}).$$

Here the implied constants depend on $d^{\#}, K$ and $M, N, \varepsilon, \delta$. But $d^{\#}$ and $K$ in turn depend only on $M, N, \varepsilon, \delta$. Thus we conclude that in (4.47), the implied constants depend only on $M, N, \varepsilon, \delta$.

In view of inequality (4.38), the right side of (4.47) is bounded below as follows:
$$(4.48)$$
$$\inf_{\{\theta \in \mathbb{R}/\mathbb{Z}:\, d\theta=0\}} H(e^{2\pi i\theta}) \geq C^{-1} \inf_{\{\theta \in \mathbb{R}/\mathbb{Z}:\, d\theta=0\}} \int_B \frac{|P(z)|^{\varepsilon}}{|Q_1(z) + e^{2\pi i\theta} Q_2(z)|^{\delta}} \, dV.$$

Inequality (4.36) is an immediate consequence of (4.47) and (4.48). The proof of the first sampling lemma is complete.

In practice, we need the first sampling lemma in the following more general form. Fix positive integers $M$ and $N_1, \ldots, N_J$, and a nondegenerate rational pair $(\varepsilon, \delta)$.

LEMMA 4.9 (second sampling lemma). *Fix nonnegative rational numbers, $\varepsilon$, $\delta$, and assume $\delta < 1$ and $(\varepsilon, \delta)$ is nondegenerate. Fix integers $M$, $N_j$, $1 \leq j \leq J$, and a constant $C > 0$. Let $s \in (0,1)$ be chosen so that the case $I = J = 1$ of Theorem 4 holds. Then there exists a positive integer $d$ with the following property. For all polynomials $P(z)$, $Q_j(z)$, $1 \leq j \leq J$, with $\deg(P) = M$, $\deg(Q_j) = N_j$ with $N_1 \geq \cdots \geq N_J$, satisfying $|||Q_1(Z) - Z^{N_1}||| < s/(2J)$ and $|||Q_j(Z)||| < s/(2J)$ for $j > 1$, we have*

$$(4.49) \qquad \int_B \frac{|P(z)|^{\varepsilon}}{(\sum_{j=1}^J |Q_j(z)|)^{\delta}} \sim \inf_{\{\theta \in (\mathbb{R}/\mathbb{Z})^J\}} \int_B \frac{|P(z)|^{\varepsilon}}{|\sum_{j=1}^J Q_j(z) e^{2\pi i\theta_j}|^{\delta}}$$
$$\sim \inf_{\{\theta \in (\mathbb{R}/\mathbb{Z})^J:\, d\theta=0\}} \int_B \frac{|P(z)|^{\varepsilon}}{|\sum_{j=1}^J Q_j(z) e^{2\pi i\theta_j}|^{\delta}}.$$

*Here all implied constants depend only on $M$, $N_j$, $\varepsilon$, $\delta$, and $s$.*

*Proof.* The case $J = 2$ is the preceding sampling Lemma 4.6. For general $J$, we observe that the right-hand side of the desired equivalence remains unchanged if we replace $(\mathbb{R}/\mathbb{Z})^J$ by $\{\theta \in (\mathbb{R}/\mathbb{Z})^J : \theta_1 = 0\}$. So we make this replacement, and then use induction. $\qquad\square$

G. *Proof of Theorem* 4. We are now in position to prove Theorem 4 in full generality. Clearly we may assume $I = 1$. As in the proof of the case $|I| = |J| = 1$, we may assume $\delta < 1$. Let $s \in (0, 1)$ be chosen as in this case (and hence as in the sampling lemmas). We may replace $Q_j(Z)$, $j \geq 2$, by $(s/(4CJ))Q_j(Z)$, without changing the size of $\sum_{j=1}^{J} |Q_j(z)|$ (up to constants depending only on $s$, $J$, and $C$). Thus we may assume that $|||Q_j(Z)||| < s/(4J)$ for $j \geq 2$. Assume that $|||Q_1(Z) - Z^{N_1}||| < s/(4J)$. Then we can apply Lemma 4.9 and write, for $\mu < s/(4J)$,

(4.50)

$$\int_B \frac{|P(z)|^\varepsilon}{(\sum_{j=1}^{J} |Q_j(z)| + \mu)^\delta} dV$$

$$\sim \inf_{\theta \in (\mathbb{R}/\mathbb{Z})^J} \int_B \frac{|P(z)|^\varepsilon}{|Q_1(z) + \sum_{j=2}^{J} Q_j(z) e^{2\pi i \theta_j} + \mu e^{2\pi i \theta_1}|^\delta} dV$$

$$\sim \inf_{\{d\theta=0\}} \int_B \frac{|P(z)|^\varepsilon}{|Q_1(z) + \sum_{j=2}^{J} Q_j(z) e^{2\pi i \theta_j} + \mu e^{2\pi i \theta_1}|^\delta} dV$$

$$\sim \inf_{\{d\theta=0\}} \int \frac{|P(z)|^\varepsilon}{|\sum_{j=1}^{J} Q_j(z) e^{2\pi i \theta_j} + \mu|^\delta} dV = \inf_{\{d\theta=0\}} \int \frac{|P(z)|^\varepsilon}{|Q_\theta(z) + \mu|^\delta} dV,$$

where $Q_\theta(z) = \sum_{\nu=0}^{N_1} b''_\nu(\theta) z^\nu$ is defined by the last equation:

$$b''(\theta) = \sum_{j=1}^{l} b''_j e^{2\pi i \theta_j}$$

(the last equation is interpreted as an identity of vectors of length $N_1$ as follows: $b''_j = (b''_{j,0}, \ldots, b''_{j,N_j})$ is identified with $(b''_{j,0}, \ldots, b''_{j,N_j}, 0, \ldots, 0)$, where $N_1 - N_j$ zeros have been added at the end to create a vector of length $N_1$).

According to Lemma 4.5, for a fixed $\theta$, the polynomial $Q_\theta(z) + \mu$ has distinct roots provided $\mu$ is sufficiently small. Thus, applying Theorem 5, we get

$$\int \frac{|P(z)|^\varepsilon}{|Q_\theta(z) + \mu|^\delta} dV \sim T(b', b''(\theta); \mu)$$

where $T$ is defined as in (4.14). Now fix $\theta \in (\mathbb{R}/\mathbb{Z})^J$ such that $d\theta = 0$. Define $B''(\theta) = \sum_{j=1}^{l} B''_j e(\theta_j)$, where we set $B''_j = (B''_{j,0}, \ldots, B''_{j,N_j}, 0, \ldots, 0)$. Define

$$\tilde{T}(B; W) = \tilde{T}(B', B''; W) = \left( \sum_{d\theta=0} T(B', B''(\theta); W)^{-1} \right)^{-1}.$$

Then, for $\mu$ sufficiently small,

$$\tilde{T}(b,\mu) \sim \inf_{d\theta=0} T(b',b''(\theta),\mu).$$

Using (4.50) we conclude that for $\mu$ sufficiently small,

$$\int \frac{|P(z)|^\varepsilon}{(\sum_{j=1}^J |Q_j(z)| + \mu)^\delta} dV \sim \tilde{T}(b,\mu).$$

Now we write $T(B;W)$ in the form (4.27), and conclude the proof exactly as we did in Section E. This completes the proof of Theorem 4.

## 5. Some preliminary considerations about stability

We now come to the application of the previous estimates to the problem of stability for inequalities of the form

$$(5.1) \qquad \int_{B^{(n)}} \frac{1}{(\sum_{j=1}^J |f_j(z)|^2)^{\delta/2}} dV_1 \cdots dV_n < \infty,$$

where the $f_j(z)$, $1 \le j \le J$, are holomorphic functions in a neighborhood of the origin in $\mathbb{C}^n$.

It is convenient to introduce

- the $J$-dimensional vector $f(z) = (f_j(z))_{j=1}^J$ and its length $|f(z)|^2 = \sum_{j=1}^J |f_j(z)|^2$,

- the polydisk $B_r^{(n)} = \{(z_1, \cdots, z_n) \in \mathbb{C}^n; |z_i| < r\}$. As we did earlier, we shall write simply $B^{(n)}$ for the polydisk of radius 1.

- the supremum norm for the $J$-vector valued functions $f(z)$

$$(5.2) \qquad ||f||^2 = \sup_{B^{(n)}} |f(z)|^2 = \sup_{B^{(n)}} \sum_{j=1}^J |f_j(z)|^2.$$

In this paper, we shall consider the following two notions of stability:

*Definition* 5.1.   (i) Inequality (5.1) is said to be uniformly stable if, for any $f_i(z)$ satisfying (5.1), and for any $0 < r < 1$, the following property holds: The functional

$$(5.3a) \qquad g = (g_j)_{j=1}^J \quad \longrightarrow \quad \int_{B_r^{(n)}} \frac{1}{(\sum_{j=1}^J |g_j(z)|^2)^{\delta/2}} dV_1 \cdots dV_n,$$

defined for holomorphic functions $g = (g_i(z))_{j=1}^J$ on $B^{(n)}$, is continuous at $f = (f_j)_{j=1}^J$ with respect to the norm (5.2).

(ii) Inequality (5.1) is said to be holomorphically stable for $d$-dimensional deformations, if for any $f(z, c) = (f_j(z, c))_{j=1}^J$ holomorphic in $B^{(n)} \times B^{(d)} \subset \mathbb{C}^n \times \mathbb{C}^d$ satisfying the inequality

$$(5.3b) \qquad \int_{B^{(n)}} \frac{1}{(\sum_{j=1}^J |f_j(z, 0)|^2)^{\delta/2}} dV_1 \cdots dV_n < \infty,$$

and for any $0 < r < 1$, the function

$$(5.4) \qquad c \quad \longrightarrow \quad \int_{B_r^{(n)}} \frac{1}{(\sum_{j=1}^J |f_j(z, c)|^2)^{\delta/2}} dV_1 \cdots dV_n$$

is continuous for $c$ in a polydisk $B_\rho^{(d)}$ in $\mathbb{C}^d$, for some $\rho > 0$.

Note that uniform stability implies holomorphic stability.

We observe that stability is not an issue when the integrals diverge, in view of Fatou's lemma. It is convenient to gather in this section the basic facts which will be required later.

A. *Order of vanishing and integrability of* $|f(z)|^{-\delta}$. Let $f(z)$ be a holomorphic function in a neighborhood of the origin in $\mathbb{C}^n$. We say that $N$ is the order of vanishing of $f(z)$ at the origin if $f(z)$ can be expressed as

$$f(z) = \sum_{i \geq N} P_i(z)$$

where for each $i$, $P_i(z)$ is a homogeneous polynomial of degree $i$, and $P_N(z)$ is not identically zero. Equivalently, the order of vanishing can be characterized as the smallest integer among the orders of vanishing at 0 of the function $f(z)$ along all lines passing through the origin. The order of vanishing $N$ of a function $f(z)$ is clearly invariant under multiplication by a nonvanishing factor.

As in the real case [18], the following simple observation is crucial to our considerations.

LEMMA 5.1.   *Assume that* $f_j(z)$, $1 \leq j \leq J$, *are holomorphic functions on a polydisk* $B_r^{(n)}$ *and that* $|f(z)|^{-\delta}$ *is integrable on* $B_r^{(n)}$ *for some* $r > 0$. *Then*

$$(5.5) \qquad\qquad \delta < \frac{2n}{N},$$

*where* $N$ *is the lowest of the orders of vanishing of all the functions* $f_j(z)$ *at the origin.*

*Proof.* Indeed, say $f_1(z)$ is a function whose order of vanishing is $N$. The fact that $N$ is the order of vanishing of $f_1(z)$ at the origin implies that $|f_1(z)| \sim |z|^N$, for $|z|$ small and $z$ in some open cone $\Gamma$ with vertex at the origin. Since all other functions $f_j(z)$, $2 \leq j \leq J$, have orders of vanishing $\geq N$, it follows that $|f_j(z)| \leq C|z|^N$ for $|z|$ small. Thus for $r' < r$ and small enough, we have $|f(z)|^2 = \sum_{j=1}^{J} |f_j(z)|^2 \sim |z|^{2N}$ for $z \in \Gamma$ and $|z| < r'$. In particular, changing to polar coordinates:

$$\infty > \int_{\Gamma \cap \{|z| < r'\}} \frac{1}{|f(z)|^\delta} dV_1 \cdots dV_n \sim \int_{\Gamma \cap \{|z|=1\}} \int_0^{r'} \rho^{2n-1-N\delta} d\rho \ ,$$

from which our assertion follows.

B. *Reduction to Weierstrass polynomials.* We say that a function $Q(z', z_n)$ is a Weierstrass polynomial in $z_n$ of degree $N$ on a polydisk $B_{r'}^{(n-1)} \times B_r^{(1)} \subset \mathbb{C}^n$ if it can be expressed as

$$(5.6) \qquad Q(z', z_n) = z_n^N + \sum_{k=0}^{N-1} a_{N-k}(z') z_n^k$$

where the coefficients $a_k(z')$, $1 \leq k \leq N$, are holomorphic on $B_{r'}^{(n-1)}$ (Note that we do no assume that $a_k(0) = 0$.) The following version of the Weierstrass Preparation Theorem provides conditions under which a function $g(z)$ can be reduced to a Weierstrass polynomial, with continuous dependence on parameters. A proof can be found in [18].

LEMMA 5.2. *Let $f(z)$ be a holomorphic function on a polydisk $B_r^{(n)}$ in $\mathbb{C}^n$. Assume that*

$$f(0, z_n) = c\, z_n^N + O(z_n^{N+1}), \ c \neq 0.$$

*Then there exists $\rho > 0$ and $0 < s < r$ with the following property. If $g(z)$ is holomorphic on $B_r^{(n)}$ and if $\sup_{B_r} |f - g| < \rho$, then $g(z)$ can be written on $B_s^{(n)}$ uniquely as*

$$(5.7) \qquad g(z) = u_g(z) Q_g(z', z_n),$$

*where $Q_g(z', z_n)$ is a Weierstrass polynomial on $B_s^{(n)}$ and $u_g(z)$ is a nonvanishing holomorphic function. Furthermore, $u_g$ and $Q_g$ depend continuously on $g$ in the following sense:*

- *For every $\tau > 0$, there exists $\tau' > 0$ such that $\sup_{B_r^{(n)}} |f - g| < \tau'$ implies that $\sup_{B_s^{(n)}} |Q_f - Q_g| < \tau$, $\sup_{B_s} |u_f - u_g| < \tau$, $\sup_{B_s^{(n)}} |1 - u_f/u_g| < \tau$, and, for all $k$, $\sup_{B_s^{(n-1)}} |a_{f,k} - a_{g,k}| < \tau$.*

- *There exist $P_0, \cdots, P_N \in \mathbb{Z}[B_{-1}, B_{-2}, \cdots][[B_1, \cdots, B_N]]$ such that the coefficients $a_{g,k}(z')$ can be expressed as a convergent sum $a_{g,k}(z') = P_k(\{b_{g,k}(z')\})$ if $g(z)$ is written in the form $g(z) = \sum_{k=0}^{\infty} b_{g,N-k}(z') z_n^k$ with $b_{g,0} = 1$. (Here $P_i$ is a power series in $B_1, \ldots, B_N$ whose coefficients are polynomials in $B_{-1}, B_{-2}, \ldots$.)*

In practice, we shall need a generalized version of Lemma 5.2, which allows us to replace expressions of the form $\sum_{j=1}^{J} |f_j(z)|^2$ by equivalent and similar expressions involving Weierstrass polynomials.

LEMMA 5.3.    Let $f(z) = (f_j(z))$, $1 \leq j \leq J$, be holomorphic functions on the polydisk $B^{(n)} \subset \mathbb{C}^n$, and let $N \geq 1$ be the lowest order of vanishing of the functions $f_i(z)$ at the origin. Then for any $s > 0$, there exists $\rho, C > 0$, and $0 < r < 1$, with the following property: After making a rotation of our coordinate system, there exists, for any $g = (g_j(z))$, $1 \leq j \leq J$, holomorphic on $B^{(n)}$ with $||g - f|| < \rho$, Weierstrass polynomials $Q_{g,j}(z', z_n)$, $1 \leq j \leq J$, of degree $N$ in $z_n$, holomorphic on the polydisk $B_r^{(n)}$, such that

$$(5.8) \qquad \sum_{j=1}^{J} |g_j(z)|^2 \quad \sim \quad \sum_{j=1}^{J} |Q_{g,j}(z)|^2$$

$$(5.9) \qquad |||Q_{g,1}(z', Z) - Z^N||| < s, \quad |||Q_{g,j}(z', Z)||| < C \quad \text{for } j \geq 2.$$

The implied constants in (5.8) depend only on $f$ and $s$, and are independent of $g$ and $z$. Furthermore, $Q_{g,i}(z)$ depends continuously on the $g$, in the sense that

- $\sup_{z' \in B_r^{(n-1)}} |||Q_{g,j} - Q_{\tilde{g},j}|||$ can be made arbitrarily small by making $||g - \tilde{g}||$ small;

- If $g = g(z, c)$ is holomorphic on the polydisk $B^{(n)} \times B^{(d)}$, then $Q_{g,j}$ is holomorphic for $(z', c)$ in a polydisk $B_\kappa^{(n+d-1)}$ for some $0 < \kappa < 1$.

*Proof.* Without loss of generality, we may assume that $N$ is the order of vanishing of $f_1(z)$ at the origin. Rotating our coordinate system, we may assume that $f_1(0, z_n) = z_n^N + O(|z_n|^{N+1})$. In view of the estimate $|f_1|^2 + |f_i|^2 \sim |f_1|^2 + |f_1 + f_i|^2$, we may replace $f_i$ by $f_1 + f_i$ if necessary, and assume that $f_i(0, z_n) = c_i z_n^N + O(|z_n|^{N+1})$, with $c_i \neq 0$. Note that all $f_i(z)$, $1 \leq i \leq I$, have now order of vanishing exactly equal to $N$ at the origin. Applying Lemma 5.2, there exists $\tilde{r} > 0$ and $\rho > 0$, so that any $g_j(z)$, $1 \leq j \leq J$, holomorphic in $B^{(n)}$ with $||g - f|| < \rho$, can be expressed as in (5.7), with $Q_{g,j}(z)$ a Weierstrass polynomial of degree $N$ on $B_{\tilde{r}}^{(n)}$.

The continuous dependence of $Q_{g,i}(z)$ on $g_i$ as stated in Lemma 5.3 is a direct consequence of the corresponding properties of Weierstrass polynomials stated in Lemma 5.2. To obtain (5.9), observe that $Q_{f,1}(0, z_n) = z_n^N$, since $f(0, z_n)$ has a zero of order $N$ at $z_n = 0$. It follows that we can make $|||Q_{f,1}(z', Z) - Z^N|||$ arbitrarily small for $z' \in B_{r'}^{(n-1)}$ by taking $0 < r'$ small enough. By the continuous dependence of the Weierstrass polynomial $Q_g(z)$ on $g$, the same is true for $|||Q_{g,1}(z', Z) - Z^N|||$ if we also take $0 < \rho$ small enough. Let $r = \min(\rho, r')$. This establishes the first part in (5.9). The second part is trivial since the Weierstrass polynomials $Q_{g_i}(z)$ have uniformly bounded coefficients, again in view of their continuous dependence on $g_i$.

Finally, (5.8) follows from the continuous dependence of $u_g$ on $g$.

C. *Uniform boundedness vs. continuity of $\int |f|^{-s}$*. In this section, we show that the continuity of the functional is, for all practical purposes, equivalent to the seemingly weaker property that

$$(5.10) \qquad \sup_{\|f - f_0\| < \rho} \int_{B_r^{(n)}} \frac{1}{|f(z)|^\delta} dV_1 \cdots dV_n < \infty,$$

for some $\rho > 0$. For our purposes, we need a somewhat more general formulation, involving quotients.

Let $r > 0$ and let $\mathcal{F} = \{(f, g)\}$ be a family of pairs of continuous functions on the polydisk $B_r^{(n)}(\mathbb{R}^n)$ in $\mathbb{R}^n$, with the property that $\{z : f(z) = 0\}$ and $\{z : g(z) = 0\}$ each have measure zero.

LEMMA 5.4. *Let $\varepsilon > 0$. Assume that there exists a constant $C$ such that*

$$(5.11) \qquad \sup_{(f,g) \in \mathcal{F}} \int_{B_r^{(n)}(\mathbb{R}^n)} \frac{|f(z)|^{1-\varepsilon}}{|g(z)|^{1+\varepsilon}} dV_1 \cdots dV_n \le C.$$

*Let $(f_0(z), g_0(z)) \in \mathcal{F}$. Then for every $\tau > 0$ there exists $\rho > 0$ so that if $(f, g) \in \mathcal{F}$ and $\sup_{B_r(\mathbb{R}^n)} |g - g_0| + |f - f_0| < \rho$, then*

$$(5.12) \qquad \Big| \int_{B_r(\mathbb{R}^n)} \frac{|f(z)|}{|g(z)|} dV_1 \cdots dV_n - \int_{B_r(\mathbb{R}^n)} \frac{|f_0(z)|}{|g_0(z)|} dV_1 \cdots dV_n \Big| < \tau.$$

*Proof.* We begin by observing that for all $(f, g)$ in $\mathcal{F}$ and all $\alpha > 0$

$$(5.13) \qquad \int_{\{|g| \le \alpha\}} |f|^{1-\varepsilon} dV_1 \cdots dV_n \le C \, \alpha^{1+\varepsilon},$$

by the Chebychev inequality. Next, we apply Hölder's inequality to the function

$$\frac{|f|}{|g|} \;=\; \Big[ \frac{|f|^{\frac{1-\varepsilon}{1+\varepsilon}}}{|g|} \Big] \cdot \Big[ |f|^{\frac{2\varepsilon}{1+\varepsilon}} \Big]$$

and exponents $p = 1 + \varepsilon$, $q = \frac{1+\varepsilon}{\varepsilon}$ to obtain

$$(5.14) \qquad \int_{\{|g| \le \alpha\}} \frac{|f|}{|g|} \le \Big( \int_{\{|g| \le \alpha\}} \frac{|f|^{1-\varepsilon}}{|g|^{1+\varepsilon}} \Big)^{\frac{1}{1+\varepsilon}} \Big( \int_{\{|g| \le \alpha\}} |f|^2 dV_1 \cdots dV_n \Big)^{\frac{\varepsilon}{1+\varepsilon}}.$$

Since we have the easy estimate

$$\int_{\{|g| \le \alpha\}} |f|^2 dV_1 \cdots dV_n \le (\sup |f|^{1+\varepsilon}) \int_{\{|g| \le \alpha\}} |f|^{1-\varepsilon} dV_1 \cdots dV_n,$$

the estimates (5.11) and (5.13) yield

$$(5.15) \qquad \int_{B_r(\mathbb{R}^n) \cap \{|g| \le \alpha\}} \frac{|f|}{|g|} dV_1 \cdots dV_n \le C' \alpha^\varepsilon,$$

with a constant $C'$ independent of $(f, g)$ in $\mathcal{F}$, if $|f|$ is uniformly bounded. Fix now $(f_0, g_0)$ in $\mathcal{F}$. Restricting ourselves to a bounded neighborhood of $(f_0, g_0)$ in the sup norm, we may assume $|f|$ is uniformly bounded. For any $\tau > 0$, we may choose and fix $\alpha$ so that the right-hand side of (5.15) is $< \frac{\tau}{3}$. Choose $\rho > 0$ so that

$$(5.16) \quad \sup_{B_r(\mathbb{R}^n)} |f - f_0| + |g - g_0| < \rho$$

$$\implies \sup_{B_r(\mathbb{R}^n) \cap \{|g_0| > \alpha/2\}} \left| \frac{|f|}{|g|} - \frac{|f_0|}{|g_0|} \right| < \frac{\tau}{3} \text{Vol}(B_r(\mathbb{R}^n))^{-1}.$$

The integral of the right-hand side of (5.12) over the region $B_r(\mathbb{R}^n) \cap \{|g_0| > \alpha/2\}$ is thus $< \tau/3$. Taking $\rho < \alpha/2$ if necessary, we can also guarantee that

$$\int_{\{|g_0| \le \alpha/2\}} \frac{|f|}{|g|} dV_1 \cdots dV_n \le \int_{\{|g| \le \alpha\}} \frac{|f|}{|g|} dV_1 \cdots dV_n,$$

for all $(f, g) \in \mathcal{F}$. The desired estimate follows.

The replacement of $\frac{|f|}{|g|}$ by $\frac{|f|^{1-\varepsilon}}{|g|^{1+\varepsilon}}$ in Lemma 5.4 is actually harmless, in view of the following lemma which extends a result of Stein [21].

LEMMA 5.5.    Let $B \subseteq \mathbb{C}^n$ be on open subset and and $B' \subseteq B$ an open relatively compact subset. Let $R$ and $S$ be functions of the form

$$R = \frac{\sum_{i=1}^{I} |f_i|^\varepsilon}{\sum_{j=1}^{J} |g_j|^\delta} \quad \text{and} \quad S = \frac{\sum_{i=1}^{I^*} |f_i^*|^{\varepsilon^*}}{\sum_{j=1}^{J^*} |g_j^*|^{\delta^*}},$$

where $f_i(z), g_j(z), f_i^*(z), g_j^*(z)$ are holomorphic functions on $B$. Assume that

$$\int_B R \ dV \quad < \quad \infty.$$

Then there exists $\sigma > 0$ such that

$$\int_{B'} R \cdot S^\sigma \ dV \quad < \quad \infty.$$

Remark. Only the fact that the functions $f_i(z), g_j(z)$ are real-analytic is actually required, but we shall not insist on this point.

Proof. Clearly we may assume that $I = I^* = 1$. Also, since

$$\int_B \frac{|f|^\varepsilon}{\sum_{j=1}^{J} |g_j|^\delta} \ dV \quad \sim \int_{[0,1]^J} \int_U \frac{|f|^\varepsilon}{|\sum_{j=1}^{J} e^{2\pi i \theta_j} g_j|^\delta} \ dV \ d\theta$$

we see, using Lemma 4.4, that there exists $\theta \in [0,1]^J$ such that

$$\int_{[0,1]^J} \frac{|f|^\varepsilon}{|\sum_{j=1}^J e^{2\pi i \theta_j} g_j|^\delta} \, dV \quad < \quad \infty.$$

Since $R < \frac{|f|^\varepsilon}{|\sum_{j=1}^J e^{2\pi i \theta_j} g_j|^\delta}$, this shows that we may assume that $J = 1$. Finally, we easily see that there exists $\theta^* \in [0,1]^{J^*}$ such that $\sum_{j=1}^{J^*} e^{2\pi i \theta_j^*} g_j^* \neq 0$. Since $S < \frac{|f^*|^{\varepsilon^*}}{|\sum_{j=1}^{J^*} e^{2\pi i \theta_j^*} g_j^*|^\delta}$ we may assume that $J^* = 1$. Thus we are reduced to proving the lemma in the case: $I = I^* = J = J^* = 1$.

Before proceeding with the rest of the proof, we recall a special case of Hironaka's theorem on resolution of singularities [9]:

HIRONAKA'S THEOREM. *Let $f$ be a holomorphic function in a neighborhood of 0 in $\mathbb{C}^n$. Then there is a complex manifold $M$ and a proper analytic map $\pi$ from $M$ to a neighborhood of 0 in $\mathbb{C}^n$ so that*

(a) *$\pi^{-1}$ is a local isomorphism from $M - \pi^{-1}(f^{-1}(0))$ to $\mathbb{C}^n - f^{-1}(0)$.*

(b) *For each $P \in M$, there are local holomorphic coordinates $z_1, \cdots, z_n$ centered at $P$ so that*

$$f \circ \pi = U \prod_{i=1}^n z_i^{a_i},$$

*where $a_i$ are nonnegative integers, and $U(z)$ is a holomorphic, nonvanishing function.*

(c) *The local coordinate system of (b) may be chosen so that*

$$J(\pi) = \tilde{U} \prod_{i=1}^n z_i^{b_i},$$

*where $J(\pi)$ is the Jacobian determinant of the map $\pi$, the $b_i$ are nonnegative integers, and $\tilde{U}(z)$ is a holomorphic, nonvanishing function.*

Parts (a) and (b) are contained in [3, Remark 8]. Part (c) follows from (a) and (b): By part (a), if $J(\pi)(z) = 0$, then $z \in \pi^{-1}(f^{-1}(0))$. Thus the hypersurface defined by $J(\pi) = 0$ is contained in the union of the coordinate hyperplanes, i.e., the hyperplanes $z_i = 0, 1 \leq i \leq n$. Thus, by the weak nullstellensatz (see [7]), the factorization of $J(\pi)$ (in the ring of germs of holomorphic functions at $z = 0$ ) must be of the form described in (c).

We return now to the proof of Lemma 5.5. For every point $p \in B$, we choose a ball $B_p \subseteq B$, centered at $p$, and apply resolution of singularities to the function $fgf^*g^*$ on $B_p$. Shrinking $B_p$ if necessary, we obtain $\pi_p : M_p \to U_p$ satisfying conditions (a), (b), (c). Let $W_p \subseteq B_p$ be the ball, centered at p, whose radius is half that of $B_p$. Then, by compactness, a finite number of the $W_p$ cover $\overline{B'}$. Thus it suffices to show that

$$\int_{W_p} R \cdot S^\sigma \; dV \quad < \quad \infty.$$

By the change of variables theorem, this amounts to showing

$$\int_{\pi_p^{-1}(W_p)} \{(R \cdot S^\sigma) \circ \pi\} \cdot |J(\pi)| \; dV \quad < \quad \infty,$$

for some $\sigma > 0$, given that the integral is finite when $\sigma = 0$. Since $\pi$ is proper, we know that $\pi_p^{-1}(\overline{W_p})$ is compact, and is thus covered by a finite number of coordinate neighborhoods which satisfy conditions (b) and (c). Thus we are reduced to proving the lemma in the case where $R$ and $S$ are of the form $\prod_{i=1}^n z_i^{c_i}$, with $c_i \in \mathbb{Z}$. In this case, the integral factors into a product of one dimensional integrals, and we are reduced to showing that if $\int_B |z|^c \; dV < \infty$, where $B \subseteq \mathbb{C}$ is a bounded open set, then $\int_B |z|^{c-\varepsilon} \; dV < \infty$ for all sufficiently small $\varepsilon \in \mathbb{R}$, which follows trivially upon changing to polar coordinates. This completes the proof of Lemma 5.5.

## 6. Stability in lower dimensions

In two complex dimensions, the stability of integrals holds in the most general form, as shown by Tian [22]. In this section, we apply our uniform one-dimensional estimates to give the following new proof of Tian's result, as well as an extension to three complex dimensions.

THEOREM 6.   Let $f_j(z)$, $1 \leq j \leq J$, be holomorphic functions on the polydisk $B^{(n)}$ of radius 1 in $\mathbb{C}^n$, and assume that (5.1) holds. Then for any $0 < r < 1$, the functional (5.3a) is continuous at $f(z) = (f_j(z))_{j=1}^J$ in the space of holomorphic functions on $B^{(n)}$, in the norm $\|f\|$ of (5.2), under any of the following assumptions:

- $n = 1$;
- $n = 2$;
- $n = 3$ and $\delta < \frac{4}{N}$, where $N$ is the lowest of the orders of vanishing of $f_j(z)$, $1 \leq j \leq J$.

*Proof.* We first assume the theorem in the case $J = 1$, and show how this implies the result for arbitrary $J$. Our hypothesis implies

$$(6.1) \qquad \int_{B^{(n)}} \frac{1}{(\sum_{j=1}^J |f_j(z)|^2)^{\delta/2}} dV_1 \cdots dV_n < \infty.$$

Lemma 5.5 implies that for every $r_1 < 1$, there exists $\delta < \delta' < 1$ such that

$$(6.2) \qquad \int_{B_{r_1}^{(n)}} \frac{1}{(\sum_{j=1}^{J} |f_j(z)|^2)^{\delta'/2}} dV_1 \cdots dV_n < \infty.$$

Lemma 4.4 implies that there exist $\theta_1, ..., \theta_J \in \mathbb{R}/\mathbb{Z}$ such that

$$(6.3) \qquad \int_{B_{r_1}^{(n)}} \frac{1}{|\sum_{j=1}^{J} f_j(z) e^{2\pi i \theta_j}|^{\delta'}} dV_1 \cdots dV_n < \infty.$$

Thus the theorem in the case $J = 1$ implies that for every $r_2 < r_1$, there exists $\rho > 0$ such that

$$(6.4) \qquad \sup_{||g-f|| < \rho} \int_{B_{r_2}^{(n)}} \frac{1}{|\sum_{j=1}^{J} g_j(z) e^{2\pi i \theta_j}|^{\delta'}} dV_1 \cdots dV_n < \infty.$$

which in turn implies

$$(6.5) \qquad \sup_{||g-f|| < \rho} \int_{B_{r_2}^{(n)}} \frac{1}{(\sum_{j=1}^{J} |g_j(z)|)^{\delta'}} dV_1 \cdots dV_n < \infty.$$

Now we apply Lemma 5.4, with $f = f_0 = 1$, $g_0 = (\sum |f_j|^2)^{\delta/2}$, $g = (\sum |g_j|^2)^{\delta/2}$ to conclude that uniform stability holds.

Now we prove Theorem 6 in the case $J = 1$. We shall write $f = f_1$.

As in Section 5.A, let $N$ be the order of vanishing of the function $f(z)$ at the origin. We saw there that the hypothesis (5.1) implies that $\delta < \frac{2n}{N}$. Thus $\delta < \frac{4}{N}$ in all three cases considered in Theorem 6. In view of Lemma 5.5, we may actually work with $\delta'$ slightly larger than $\delta$ but still satisfying $\delta' < \frac{4}{N}$. It suffices then, in view of Lemma 5.4 and the observation subsequent to Lemma 5.5, to show that the integrals of $|g(z)|^{-\delta'}$ are uniformly bounded for $||g - f||$ small enough. Let $s \in (0, 1)$ be as in Theorem 4, and let $r > 0$, $\rho > 0$ be chosen as in Lemma 5.2. In particular, to each $g(z)$ corresponds a Weierstrass polynomial $Q_g(z)$ of degree $N$ in $z_n$, and

$$(6.6) \qquad \int_{B_r^{(n)}} \frac{1}{|g(z)|^{\delta'}} dV_1 \cdots dV_n \sim \int_{B_r^{(n)}} \frac{1}{|Q_g(z)|^{\delta'}} dV_1 \cdots dV_n.$$

Changing scales, we may assume that $B_r^{(n)}$ is the polydisk $B^{(n)}$ of radius 1.

Now $\delta' < \frac{4}{N}$ and all the hypotheses in Theorem 3, (b), are verified. We deduce that

$$(6.7) \qquad \int_{B^{(1)}} \frac{1}{|Q_g(z)|^{\delta'}} dV_n \sim \frac{1}{(\sum_{j=1}^{J^*} |G_{g,j}(z')|)^{\delta'}},$$

where the expressions $G_{g,j}(z')$ are polynomials in the coefficients of the Weierstrass polynomial $Q_g(z', Z)$, and the index $g$ indicates their dependence on the function $g$. In particular, $G_{g,j}(z')$ depends continuously on $g(z)$. We obtain:

$$(6.8) \qquad \int_{B^{(n)}} \frac{1}{|g(z)|^{\delta'}} dV_1 \cdots dV_n \sim \int_{B^{(n-1)}} \frac{1}{(\sum_{j=1}^{J^*} |G_{g,j}(z')|)^{\delta'}} dV_1 \cdots dV_{n-1}.$$

We consider successively the three cases listed in Theorem 6.

*Case* 1.   In the case $n = 1$, no additional integration besides $dV_1$ is necessary, and we just consider (6.7). Since $|f(z)|^{-\delta'}$ is integrable, the right-hand side of (6.3) is finite for $g(z) = f(z)$. This means that $G_{f,j} \neq 0$ for some $j$. Clearly $|G_{g,j}|$ is then bounded away from 0 for $\|g - f\|$ small enough. Applying (6.7) again gives the theorem in this case.

*Case* 2.   In the case $n = 2$, the right-hand side of (6.8) is an integral of the same form as the original integral, but in dimension 1. Since this case is now known to be stable, the uniform boundedness of the integrals on the left-hand side of (6.8) follows.

*Case* 3.   The argument is identical to Case 2, with the uniform boundedness of the 3-dimensional integrals on the left-hand side of (6.8) reduced to that of the two-dimensional integrals on the right-hand side, which has just been established.                                                                 □

## 7. Holomorphic stability in arbitrary dimensions

We turn now to the stability of integrals of the form (5.1) in arbitrary dimensions. Very little is known in this case. There are no counterexamples such as Varchenko's counterexample in 3 real dimensions. The only general results available so far are the following two theorems due respectively to B. Lichtin [14] and Y.-T. Siu [19], [20], which deal both with one-dimensional parameter families of holomorphic functions.

LICHTIN'S THEOREM.    *Let $f(z, c)$ be a one-parameter holomorphic family of germs of holomorphic functions at the origin in $\mathbb{C}^n$. Assume that for each $c$, the germ $f(z, c)$ has an isolated singularity at the origin, and that $\int_{B^{(n)}(0)} |f(z, 0)|^{-\delta} dV_1 \cdots dV_n < \infty$ for a Milnor ball $B^{(n)}(0)$ for $f(z, 0)$. Then for any family of Milnor balls $B^{(n)}(c) \subset\subset B^{(n)}(c)'$ for each germ $f(z, c)$, the integral $\int_{B^{(n)}(c)} |f(z, c)|^{-\delta} dV_1 \cdots dV_n$ is finite for each $c$ small enough.*

SIU'S SEMICONTINUITY LEMMA.    *Let $B^{(n)} \times B^{(1)}$ be an open polydisk centered at the origin in $\mathbb{C}^n \times \mathbb{C}$. Let $B_r^{(n)}$, $0 < r < 1$ be any relatively compact open polydisk contained in $B^{(n)}$. Let $f(z, c) = (f_j(z))$, $1 \leq j \leq J$, be holomorphic functions on $B^{(n)} \times B^{(1)}$ and $\delta$ be a positive number. If* $\int_{B^{(n)}} |f(z,0)|^{-\delta} dV_1 \cdots dV_n$ *is finite, then there exists a constant $C < \infty$ so that*

$$(7.1) \qquad \inf_{0 < |c| < \rho} \int_{B_r^{(n)}} |f(z,c)|^{-\delta} dV_1 \cdots dV_n \leq C,$$

*for all $\rho > 0$ sufficiently small.*

Equivalently, there is a subsequence $c_\nu \to 0$ so that

$$\sup_\nu \int_{B_r^{(n)}} |f(z,c_\nu)|^{-\delta} dV_1 \cdots dV_n \leq C.$$

As Siu noted [20], a more careful scrutiny of his arguments actually shows that

$$\lim_{\nu \to \infty} \int_{B_r^{(n)}} |f(z,c_\nu)|^{-\delta} dV_1 \cdots dV_n = \int_{B_r^{(n)}} |f(z,0)|^{-\delta} dV_1 \cdots dV_n.$$

Alternatively, this slight strengthening also follows automatically in view of the discussion in Section 5.C.

Siu's semicontinuity lemma is a direct consequence of a more general Lemma (henceforth referred to as Siu's lemma), which Siu establishes for nonpositive plurisubharmonic functions, using the work of Ohsawa-Takegoshi [16] on extensions of holomorphic functions with Carleman weights. For greater clarity, we present both the statement and proof of Siu's lemma in a separate appendix.

The purpose of this section is to show how the uniform estimates for ARP's derived earlier, combined with Siu's lemma on plurisubharmonic functions, can give holomorphic stability for 1-parameter deformations. In effect, the infimum on the left-hand side of (7.1) can be replaced by a supremum. This strengthening implies Lichtin's theorem, and shows that no example such as Varchenko's can exist in the complex case.

THEOREM 7. *Let $D \subseteq \mathbb{C}^n$ be an open set and $D' \subseteq D$ a relatively compact subset. Let $R(z,c)$ be a function of the form*

$$(7.3) \qquad R(z,c) = R(z,c;\varepsilon,\delta) = \frac{(\sum_{i=1}^I |f_i(z,c)|^2)^{\varepsilon/2}}{(\sum_{j=1}^J |g_j(z,c)|^2)^{\delta/2}},$$

*where $g_i(z,c)$, $f_j(z,c)$ are holomorphic functions on $D \times B^{(1)}$, where $B^{(1)}$ is a ball centered at the origin in the $\mathbb{C}$ plane, and $\varepsilon, \delta$ are fixed nonnegative rational numbers. Assume that*

$$(7.4) \qquad \int_D R(z,0) dV_1 \cdots dV_n < \infty.$$

*Assume as well that*

$(7.5)$ $-\log(R(z,c))$ *is a nonpositive plurisubharmonic function on $D \times B^{(1)}$ .*

*Then there exists a disk $B_\rho^{(1)}$ so that the function*

$$(7.6) \qquad\qquad c \quad \longrightarrow \quad \int_{D'} R(z;c) \; dV$$

*is finite and continuous on $B_\rho^{(1)}$.*

The Main Theorem in the introduction follows by taking $\sum_{i=1}^{I} |f_i(z,c)|^2$ to be a large constant.

Before giving the detailed proof of Theorem 7, we pause to discuss a global strategy for such problems, and point out some key differences between the present situation and the lower dimensional cases treated in Section 6.

By compactness, we may assume that $D$ is the polydisk $B^{(n)}$ of radius 1, and that $D'$ is another small polydisk $B_r^{(n)}$ with $r > 0$ suitably chosen. In view of our earlier observations in Section 5.C, it again suffices to prove that

$$(7.7) \qquad\qquad \sup_{c \in B_\rho^{(1)}} \int_{B_r^{(n)}} R(z,c; \varepsilon - \nu, \delta + \nu) \; dV_1 \cdots dV_n < \infty$$

for some $\rho, \nu > 0$. Let $N$ be the lowest order of vanishing of $g_j(z,0)$ at $z = 0$, and choose $s \in (0,1)$ as in Theorem 4. Applying Lemma 5.3, we can find $r > 0$, $\rho > 0$, so that $g_j(z,c)$ can be written in the polydisk $B_r^{(n)} \times B_\rho^{(1)}$ as $g_j(z,c) = u_j(z,c)Q_j(z,c)$, where $Q_j(z,c)$ is a Weierstrass polynomial in $z_n$ on $B_r^{(n)}$, depending holomorphically on $c \in B_\rho^{(1)}$, $u_j(z,c)$ is bounded away from 0, and (5.8) and (5.9) hold. The $f_i(z,c)$ may be assumed to have a similar Weierstrass factorization: $f_i(z,c) = v_i(z,c)P_i(z,c)$. It clearly suffices to establish the estimate (7.3) with $f_i$ replaced by $P_i$ and $g_j$ replaced by $Q_j$.

To simplify the notation, it is also convenient to dilate the $z$ variables, and assume that $B_r^{(n)}$ is just the unit polydisk $B^{(n)}$.

The function $\frac{(\sum_{i=1}^{I} |P_i(z)|^2)^{\varepsilon/2}}{(\sum_{j=1}^{J} |Q_j(z,c)|^2)^{\delta/2}}$ is thus of the type associated to the space $\mathcal{B}$ of (4.4). Recall that $\mathcal{B}$ is essentially a space of polynomials $(P_i(Z), Q_j(Z))$, whose coordinates are the coefficients $b'$ and $b''$ of $P_i(Z)$ and $Q_j(Z)$. Set $b = (b', b'')$ as before. Let $\mathcal{U}_\lambda$, $T_\lambda(b) \in \mathbb{Z}|b|$ be the filtration following from Theorem 4, and let $b(z', c)$ be the coefficients of the polynomials $(P_i(z', Z, c), Q_j(z', Z, c)$. Given $(z', c)$, there exists a unique $\lambda(z', c)$ so that $b(z', c) \in \mathcal{U}_{\lambda(z',c)} \setminus \mathcal{U}_{\lambda(z',c)+1}$. Theorem 4 provides then the following estimate

$$\int_{B^{(n)}} R(z,c)dV_n \sim \int_{B^{(n)}} \frac{(\sum_{i=1}^{I} |P_i(z,c)|^2)^{\varepsilon/2}}{(\sum_{i=1}^{J} |Q_i(z,c)|)^{\delta/2}} dV_n \sim T_{\lambda(z',c)}(b(z',c)).$$

Our main task is to carry out the remaining $dV_1 \cdots dV_{n-1}$ integrals. The main difficulties are the following.

• Let $\lambda$ be the stratification level characterized by the fact that $\mathcal{U}_\lambda$ is the smallest subvariety which contains $(P(z',c),Q(z',c)) = (P_i(z',c),Q_j(z',c))$ for all $(z',c)$. In general, we have then

$$T_{\lambda(z',c)}(b(z',c)) = T_\lambda(b(z',c))$$

except possibly when $(P(z',c),Q(z',c)) \in \mathcal{U}_{\lambda+1}$. The $c$-parameter space gets partioned into the Zariski open subspace $W_0$ of values of $c$'s for which the variety

$$Z_c = \{z' \in B^{(n-1)} : (P(z',c),Q(z',c)) \in \mathcal{U}_{\lambda+1}\}$$

is of codimension $\geq 1$ in $B^{(n-1)}$, and the closed variety $W_1$ of values of $c$ for which the variety $Z_c$ is of codimension 0. In the first case, the variety $Z_c$ is of measure 0, and for the purpose of integrating in $dV_1 \cdots dV_{n-1}$, we may just use $T_\lambda(z',c)$. In the second case, we must use instead $T_{\lambda+1}(z',c)$ on some open subset of $W_1 \backslash W_2 \subseteq W_1$. Continuing in this fashion, we get a stratification $W_0 \supseteq W_1 \supseteq W_2 \cdots$ of $c$-parameter space, with different $T_\lambda$ controlling the sizes of the integrals at each level of the stratification. The fact that, at intermediate stages, there is no single expression controlling the integral, but rather many expressions depending on which stratification variety we are in, is a source of significant difficulties.

• As a rule, the expression $T_\lambda(z',c)$ can be expected to get increasingly difficult to control as $c$ approaches a lower dimensional variety of the $c$-stratification. Indeed, the emergence of $T_{\lambda+1}(z',c)$ is due to the fact that $T_\lambda(z',c)$ approaches an indefinite expression $0/0$.

• It is not surprising that difficulties should surface at a second integration. The first integration produces in general a rational expression. Integrals of *generic* rational expressions are manifestly not stable under *arbitrary* deformations.

• A successful approach must identify special features of one of the following:

　　1. The rational expression appearing at intermediate stages.

　　2. The deformations resulting from deformations of the integrand $R(z,0)$.

• The lower-dimensional cases treated in Section 6 can be viewed as an example of the first strategy (cf. (6.4)). The proof of Theorem 6 below can be viewed as an example of the second, with Siu's semicontinuity theorem controlling indirectly the relevant deformations.

*Proof of Theorem* 7. Our goal is to prove the estimate (7.7) for some $\rho, \nu > 0$. We first choose $\nu$ so that if we replace $R(z,c;\varepsilon,\delta)$ by $R(z,c;\varepsilon-\nu,\delta+\nu)$, then (7.4) holds, and the pair of exponents $(\varepsilon-\nu,\delta+\nu)$ is nondegenerate.

This is possible since, by Lemma 5.5, (7.4) will hold for all sufficiently small $\nu$, and, by the definition of nondegeneracy, the pair $(\varepsilon - \nu, \delta + \nu)$ is degenerate for only a finite number of values of $\nu$. Thus we fix a $\nu$ for which these conditions are met, and we replace $R(z, c; \varepsilon, \delta)$ by $R(z, c; \varepsilon - \nu, \delta + \nu)$. We note that (7.5) will hold for this choice of $\nu$ as well (see [10, Corollaries 1.66 and 1.68]).

Next, we see by Lemma 5.3, that for the purposes of proving (7.7), we may assume, after making a linear change in the $z$ variables, that the $f_i(z', z_n; c)$ and the $g_j(z', z_n; c)$ are Weierstrass polynomials in $z_n$, for $(z; c) = (z', z_n; c)$ in the polydisk $B_{r_{n-1}}^{(n-1)} \times B_{r_n}^{(1)} \times B_\rho^{(1)}$, for some $r_{n-1}, r_n, \rho > 0$. Decreasing the size of $r_{n-1}$ and $\rho$ will guarantee that the hypothesis of Theorem 4 applies to the one dimensional integral in the $z_n$ variable over the ball $B_{r_n}^{(1)}$. Thus we obtain:

$$\int_{B_{r_n}^{(1)}} R(z; c) \ dV_n \sim T_{\lambda(b(z'; c))}(b(z'; c)) \ = \ T_{\lambda(z'; c)}(z'; c) \ .$$

Let $\lambda' = \inf_{(z'; c)} \lambda(z'; c)$ where the infimum is taken over the polydisk $B_{r_{n-1}}^{(n-1)} \times B_\rho^{(1)}$, and write $T_{\lambda'}(z'; c) = K_{\lambda'}(z'; c)/L_{\lambda'}(z'; c)$. We must have that $L_{\lambda'}(z'; c) \not\equiv 0$, for otherwise, the integral given in (7.6) would be infinite for all sufficiently small nonzero $c$, which, by Siu's lemma, would contradict (7.4). Now, for each nonzero $c \in B_\rho^{(1)}$, Theorem 4 implies that for all $z'$ outside the subvariety $Z_c' = \{z' : K_{\lambda'}(z'; c) = 0\}$:

$$(7.8) \qquad \int_{B_{r_n}^{(1)}} R(z; c) \ dV_n \sim T_{\lambda'}(z'; c) = |c|^\mu \cdot \tilde{T}_{\lambda'}(z'; c)$$

where $\mu$ is a rational number, $\tilde{T}_{\lambda'} = (\tilde{K}')^{\varepsilon'}/(\tilde{L}')^{\delta'}$, $\tilde{K}'$ and $\tilde{L}'$ are sums of absolute values of holomorphic functions, and $\tilde{K}'(z'; 0) \not\equiv 0$, $\tilde{L}'(z'; 0) \not\equiv 0$ (to see this, we simply factor out the highest power of $|c|$ from the numerator and denominator of $T_{\lambda'}$).

Now choose $a'$ such that $\tilde{K}(a'; 0) \neq 0$. Then $\tilde{K}(a'; c) \neq 0$ for $c$ sufficiently small. Hence, shrinking $\rho$ is necessary, we may assume that $Z_c'$ is a proper subvariety for all nonzero $c \in B_\rho^{(1)}$. In particular, $Z_c'$ has measure zero.

After making a linear change in the $z'$ variables if necessary, we may assume that $\tilde{K}'(z'; c)$ and $\tilde{L}'(z'; c)$ do not vanish identically when $(z'', c) = (z_1, \ldots, z_{n-2}; c) = 0$. Using Weierstrass preparation again, we may assume that $\tilde{K}'(z'; c)$ and $\tilde{L}'(z'; c)$ are sums of absolute values of monic polynomials in $z_{(n-1)}$ whose coefficients are analytic functions of $(c; z'')$ where $z'' \in B_{r_{n-2}}^{(n-2)}$. Choose $\sigma' > 0$ to be a small rational number with the following properties:

a) $(\tilde{K}')^{\varepsilon'}/(\tilde{L}')^{\delta' + \sigma'}$ is nondegenerate with respect to the $z_{n-1}$ variable.

b)

$$\int_{B_{r_{n-2}}^{(n-2)} \times B_{r_{n-1}}^{(1)} \times B_{r_n}^{(1)}} \frac{R(z; 0)}{(\tilde{L}'(z'; 0))^{\sigma'}} \ dV_1 \cdots dV_n \ < \ \infty.$$

From Lemma 5.5 we know that b) holds for all sufficiently small $\sigma' > 0$. Moreover, a) will be satisfied for all but finitely many $\sigma'$. Thus there does exist a rational $\sigma' > 0$ for which a) and b) are simultaneously satisfied.

Thus, shrinking $\rho$ and $r_2$ if necessary, for all nonzero $c \in B_\rho^{(1)}$ we have

$$(7.9) \quad \int_{B_{r_{n-1}}^{(1)}} \int_{B_{r_n}^{(1)}} \frac{R(z;c)}{(\tilde{L}'(z';c))^{\sigma'}} \, dV_n dV_{n-1} \sim |c|^\mu \cdot \int_{B_{r_{n-1}}^{(1)}} \frac{\tilde{T}_{\lambda'}(z';c)}{(\tilde{L}'(z';c))^{\sigma'}} dV_{n-1}.$$

We may apply Theorem 4 again (since the nondegeneracy hypothesis is now valid), and we find that

$$\int_{B_{r_{n-1}}^{(1)}} \frac{\tilde{T}_{\lambda'}(z';c)}{(\tilde{L}'(z';c))^{\sigma'}} dV_{n-1} \ \sim \ T_{\lambda''}(z'';c) \ \sim |c|^{\mu'} \, \tilde{T}_{\lambda''}(z'';c)$$

for $c$ sufficiently small, and for all $z'' \in B_{r_{n-2}}^{(n-2)}$ lying outside a proper closed subvariety which depends on $c$. Continuing as before, we write $\tilde{T}_{\lambda''} = (\tilde{K}'')^{\varepsilon''}/(\tilde{L}'')^{\delta''}$ and choose $\sigma'' > 0$ such that

a) $(\tilde{K}'')^{\varepsilon''}/(\tilde{L}'')^{\delta'+\sigma''}$ is nondegenerate with respect to the $z_{n-2}$ variable.

b)
$$\int_{B_{r_{n-2}}^{(n-2)} \times B_{r_{n-1}}^{(1)} \times B_{r_n}^{(1)}} \frac{R(z;0)}{(\tilde{L}'(z';0))^{\sigma'} (\tilde{L}''(z'';0))^{\sigma''}} \, dV_1 \cdots dV_n \ < \ \infty.$$

Continuing by induction, and changing notation slightly, we obtain:

$$(7.10) \qquad\qquad \int_P \frac{R(z;c)}{\Lambda(z;c)} \, dV_1 \cdots dV_n \ \sim T(c)$$

where the left-hand side of (7.10) is finite at $c = 0$ and the right-hand side of (7.10) is an ARP in one variable. Here $\Lambda(z;c) = (L'(z';c))^{\sigma'} (L''(z'';c))^{\sigma''} \cdots$, and $P$ is the polydisk $\prod_{k=1}^n B_{r_k}^{(1)}$.

One easily sees that $-\log(\Lambda(z,c))$ is a nonpositive plurisubharmonic function (see [10, Corollaries 1.66 and 1.68]). Thus we can apply Siu's lemma, which implies that the left-hand side of (7.10) is bounded on some subsequence approaching zero. In particular, $\lim_{c \to 0} T(c)$ is finite, which implies that

$$\sup_{c \in B_\rho^{(1)}} \int_P \frac{R(z;c)}{\Lambda(z;c)} \, dV_1 \cdots dV_n < \infty;$$

hence

$$\sup_{c \in B_\rho^{(1)}} \int_P R(z;c) \, dV_1 \cdots dV_n < \infty.$$

This gives estimate (7.7). Theorem 7 is proved.

## 8. Uniform estimates for distribution functions

It may be worth stating separately the following theorem, which is an easy consequence of the previous developments. Let $f(z) = (f_j(z))_{j=1}^J$ be a $J$-vector of analytic function on the unit ball $B^{(n)}(\Omega)$ in either $\Omega = \mathbb{R}^n$ or $\Omega = \mathbb{C}^n$ (in which case analyticity is the same as holomorphicity). Let

$$(8.1) \qquad \mu_f(\alpha, r) = \mathrm{Vol}\{\alpha \in B_r^{(n)}(\Omega); |f(z)| < \alpha\}.$$

THEOREM 8.   *Let $0 < r' < r < 1$. The following statements hold*:

$$(8.2) \qquad \int_{B^{(n)}(\Omega)} |f(z)|^{-\delta} dV_1 \cdots dV_n < \infty \Rightarrow \sup_{g \in \mathcal{F}} \mu_g(\alpha, r) \le C\alpha^\delta,$$

$$(8.3) \qquad \mu_f(\alpha, r) \le C_r \alpha^\delta \Rightarrow \sup_{g \in \mathcal{F}} \mu_g(\alpha, r') \le C'_{\delta^\#, r'} \alpha^{\delta^\#}, \quad \text{for all } 0 < \delta^\# < \delta,$$

*when*

(a) $\Omega = \mathbb{C}^2$ *(respectively $\mathbb{R}^2$), $\mathcal{F}$ is the family of functions $J$-vectors $g(z) = (g_j(z))_{j=1}^J$ holomorphic (respectively analytic) in $B^{(2)}(\Omega)$, and satisfying $\sup_{B^{(2)}(\Omega)} |g - f| < \rho$, for $\rho$ small enough.*

(b) $\Omega = \mathbb{C}^n$, *$\mathcal{F}$ is the family of $J$-vectors of the form $g = (g_j(z, c))_{j=1}^J$, where $g_j(z, c)$ is holomorphic in a polydisk in $\mathbb{C}^n \times \mathbb{C}$, and $|c|$ is small enough.*

(c) $\Omega = \mathbb{C}^3$ *(respectively $\mathbb{R}^3$), $\delta < \frac{4}{N}$, where $N$ is the lowest order of vanishing of the function $f_j(z)$ in the unit polydisk, $\mathcal{F}$ is the family of $J$-vectors $g(z) = (g_j(z))_{j=1}^J$ holomorphic (respectively analytic) in $B^{(3)}(\Omega)$, and satisfying $\sup_{B^{(3)}(\Omega)} |g - f| < \rho$, for $\rho$ small enough.*

*Proof.* Assume the inequality on the left of (8.2). Under any of the three sets of conditions (a)–(c), the integrals $\sup_g \int_{B^{(n)}(\Omega)} |g|^{-\delta'} dV$ are bounded ((c) and the complex cases of (a) and (b) have been established in this paper. The real cases of (a) and (c) have been established in [18].) The desired statement follows from the Chebychev inequality (1.2). Assume now the inequality on the left of (8.3). The identity

$$\int_{B_r^{(n)}(\Omega)} |f|^{-\delta^\#} dV = \frac{1}{\delta^\#} \int_{B_r^{(n)}(\Omega)} \Big( \int_{|f|}^\infty \alpha^{-\delta^\# - 1} d\alpha \Big) dV$$
$$= \frac{1}{\delta^\#} \int_0^\infty \alpha^{-\delta^\# - 1} \mu_f(\alpha, r) d\alpha$$

shows that the integral on the left-hand side is finite for all $0 < \delta^\# < \delta$. We may now apply the previous result.   $\square$

## Appendix: Siu's lemma

As we had mentioned earlier, the key lemma of Siu which we need is in an unpublished paper of his on Kähler-Einstein metrics for Fano manifolds. Since Siu's paper on Kähler-Einstein metrics is unpublished, we reproduce here with his kind permission the statement and proof of his lemma, which he communicated to us several years ago.

SIU'S LEMMA. *Let $D$ be a bounded Stein open set in $\mathbb{C}^n$. Let $D'$ be a relatively compact open subset of $D$. Let $\varphi(z,c)$ be a nonpositive plurisubharmonic function on $D \times B^{(1)}$ such that $\int_D \exp(-\varphi(z,0))dV_1 \cdots dV_n < \infty$. Then there exists a positive number $C$ depending only on $D$ and $D'$ such that*

$$(7.2) \quad \inf_{0 < |c| < \rho} \int_{D'} \exp(-\varphi(z,c))dV_1 \cdots dV_n \le C \int_D \exp(-\varphi(z,0))dV_1 \cdots dV_n$$

*for $0 < \rho$ sufficiently small.*

Siu's lemma implies the earlier Siu semicontinuity lemma described in Section 7, by setting $\varphi(z,c) = \delta \log \left( \sum_{j=1}^J |f_j(z,c)|^2 \right) - A$, where $A$ is a large positive number.

The proof of Siu's lemma is based on the following result due to Ohsawa and Takegoshi [16]:

OHSAWA-TAKEGOSHI THEOREM. *Let $\Omega$ be a bounded pseudoconvex domain in $\mathbb{C}^n$, $\psi : \Omega \to \mathbb{R} \cup \{-\infty\}$ be a plurisubharmonic function, and $H \subset \mathbb{C}^n$ be a complex hyperplane. Then there exists a constant $C$ depending only on the diameter of $\Omega$ such that, for any holomorphic function $f$ on $\Omega \cap H$ satisfying*

$$\int_{\Omega \cap H} e^{-\psi}|f|^2 dV_1 \cdots dV_{n-1} < \infty,$$

*there exists a holomorphic function $F$ on $\Omega$ satisfying $F|_{\Omega \cap H} = f$ and*

$$\int_{\Omega} e^{-\psi}|F|^2 dV_1 \cdots dV_n \le C \int_{\Omega \cap H} e^{-\psi}|f|^2 dV_1 \cdots dV_{n-1}.$$

The original proof of the Ohsawa-Takegoshi theorem is in [16]. A simpler proof due to Siu can be found in [19].

*Proof of Siu's lemma.* Let $\varphi_\eta(w,c) = \varphi(w, \eta c)$ for $0 < \eta < 1$. Apply the Ohsawa-Takegoshi result to the case $\Omega = D \times B^{(1)}$ and $f \equiv 1$ on $D \times 0$ and to the plurisubharmonic function $\varphi_\eta$. We get a holomorphic function $F_\eta(w,c)$ on $D \times B^{(1)}$ such that $F_\eta(w,0) \equiv 1$ on $D \times 0$ and

$$(A.1) \qquad \int_{D \times B^{(1)}} |F_\eta|^2 \exp(-\varphi_\eta) \le C \int_{w \in D} \exp(-\varphi(w,0))$$

with $C$ independent of $\eta$. Since $\varphi_\eta$ is nonpositive, we have

$$\int_{D\times B^{(1)}} |F_\eta|^2 \leq C \int_{w\in D} \exp(-\varphi(w,0)).$$

Let $B_\rho^{(1)}$ be the open disk $\{|c|<\rho\}$ in $\mathbb{C}$ of radius $\rho$ centered at 0. Let $D''$ be a relatively compact open neighborhood of the closure of $D'$ in $D$. Since the family of holomorphic functions $\{F_\eta\}$ is uniformly $L^2$ on $D\times B^{(1)}$, it follows that the supremum norm of $F_\eta$ and its first-order partial derivatives are uniformly bounded on $D''\times B_{1/2}^{(1)}$. Since $F_\eta(w,0)\equiv 1$ on $D\times 0$, it follows that there exists some positive number $\rho<\frac{1}{2}$ and a positive number $a$ independent of $\eta$ such that the infimum of $|F_\eta|$ on $D'\times B_\rho^{(1)}$ is at least $a$.

We rescale the $c$ variable in (A.1) by substituting $c$ for $\eta c$ and get

$$\frac{1}{\eta^2} \int_{D\times B_\eta^{(1)}} \exp(-\varphi) \leq C\, a^{-2} \int_{w\in D} \exp(-\varphi(w,0))$$

for $0<\eta<\rho$, and

$$\inf_{0<|c|<\rho} \int_{w\in D'} \exp(-\varphi(w,c)) \leq C\, a^{-2}\pi^{-1} \int_{w\in D} \exp(-\varphi(w,0))$$

for $0<\eta<\rho$.                                                                    $\square$

COLUMBIA UNIVERSITY, NEW YORK, NY
*E-mail address:* phong@math.columbia.edu

RUTGERS UNIVERSITY, NEWARK, NJ
*E-mail address:* sturm@andromeda.rutgers.edu